\documentclass[a4paper,12pt]{article}

\usepackage[T1]{fontenc}
\usepackage[utf8]{inputenc}
\usepackage{amssymb,amsmath,amsfonts,mathrsfs}
\usepackage{latexsym}
\usepackage[pagebackref,colorlinks=true,urlcolor=blue,linkcolor=red,citecolor=red]{hyperref}
\hypersetup{breaklinks=true}
\usepackage[usenames,dvipsnames]{color}
\usepackage[figuresright]{rotating}
\usepackage{times}
\usepackage{lmodern}
\usepackage{makeidx}
\usepackage[intoc]{nomencl}
\usepackage{overpic,caption,subcaption}
\usepackage{textcomp}
\usepackage{nicefrac}

\usepackage{fancyhdr}  

\usepackage[left=2cm,top=3.5cm,right=2cm,bottom=3.5cm,footskip=70pt, footnotesep=24pt]{geometry}
\usepackage{enumitem}

\usepackage{pgfpages}
\usepackage{pgfplots}
\usepackage{tikz-3dplot} 
\pgfplotsset{compat=newest}
\usetikzlibrary[patterns]
\usetikzlibrary{arrows,positioning,shapes}  
\usetikzlibrary{plotmarks}
\usetikzlibrary{decorations.pathmorphing}
\usetikzlibrary{calc}
\pgfplotsset{compat=newest}

\tikzset{
    myarrow/.style={
        draw,
        fill=red,
        single arrow,
        minimum height=4.5ex,
        single arrow head extend=1ex
    }
}

\usepackage{empheq}
\usepackage[linesnumbered,ruled,vlined]{algorithm2e}
\SetKwInput{KwInput}{Input}  
\SetKwInput{KwOutput}{Output}
\SetKw{KwBy}{by} 
\DeclareMathOperator*{\argmin}{arg\,min}

\newcommand{\bs}{\boldsymbol}

\newcommand{\wtbF}[1]{\overline{\widetilde{\bsF}}}

\newcommand{\bsF}{\boldsymbol{F}}

\newcommand{\bsx}{\boldsymbol{x}}

\usepackage{xcolor}
\usepackage{listings}
\lstdefinestyle{terminal}
{
    backgroundcolor=\color{black},
    basicstyle=\scriptsize\color{white}\ttfamily
}

\title{$r-$Adaptive Deep Learning Method for Solving Partial Differential Equations} 
\author{Ángel J. Omella$^1$ \and David Pardo$^{1,2,3}$}
\date{%
    $^1$University of the Basque Country (UPV/EHU), Leioa, Spain\\%
    $^2$Basque Center for Applied Mathematics (BCAM), Bilbao, Spain\\%
    $^3$Ikerbasque, Bilbao, Spain\\[2ex]%
    \today
}

\begin{document}

\pgfkeys{/pgf/number format/.cd,1000 sep={\,}}

\maketitle

\begin{abstract}
We introduce an $r-$adaptive algorithm to solve Partial Differential Equations using a Deep Neural Network.
The proposed method restricts to tensor product meshes and optimizes 
the boundary node locations in one dimension, from which we build two- or three-dimensional meshes.
The method allows the definition of fixed interfaces to design conforming meshes, and enables changes in the topology, i.e., some nodes can jump across fixed interfaces.
The method simultaneously optimizes the node locations and the PDE solution values over the resulting mesh.
To numerically illustrate the performance of our proposed $r-$adaptive
method, we apply it in combination with a collocation method, a Least Squares Method, and a Deep Ritz Method. 
We focus on the latter to solve one- and two-dimensional problems whose solutions are smooth, singular, and/or exhibit strong gradients.
\end{abstract}
\section{Introduction}\label{Intro}

Deep Learning (DL) \cite{bookDL,Koll2021} is nowadays applied to multiple fields \cite{DLapplications}, including biomedical applications \cite{biomedicine}, structural health monitoring \cite{puente_ana,Ana_puentes_FEM}, and geosteering \cite{extra_deep}.
Indeed, DL can perform complex tasks with high accuracy without incurring prohibitive computational costs. 
DL has allowed an essential advance in solving problems where the relationship between input and output data is complex and unknown. 
For example, merging DL techniques with the Finite Element Method can be used to improve the solution of Partial Differential Equations (PDEs) \cite{BREVIS,Brevis2,Maciej1,Maciej2,carlos}.
In addition, the use of DL to predict PDEs behavior has also raised great interest during the last decade\cite{overview,PIML,PINN}.

To solve a PDE using DL, we define a loss function whose global minimum satisfies the PDE and the boundary conditions (BCs). 
The selection of the numerical method to solve the PDE, formulated in strong, weak, or ultra-weak form, leads to different definitions of the loss function.
Some existing methods used in this context are: the Deep Ritz method (DRM) \cite{DeepRitz} based on the minimization of the energy of the PDE solution, 
Physics-Informed Neural Networks (PINNs) \cite{PINN} based on collocation methods,  
the Deep first-order Least-Squares method (DLS) \cite{DeepLS},
the Deep Galerkin method (DGM) \cite{DGM},
and the $hp$-Variational Physical Informed Neural Networks ($hp$-VPINNs) ~\cite{hpVPINNs} based on a Petrov-Galerkin domain decomposition.

The potential of Deep Neural Networks (DNNs) to solve high-dimensional PDEs has been proved, for example, in \cite{high_dimensions2, high_dimensions}. 
For these applications, Monte Carlo methods \cite{Montecarlo_no_lineal} are used for integration.
However, to deal with low-dimensional PDEs, in which the efficiency of Monte Carlo methods is limited compared to other methods, we can resort to  standard mesh-based quadrature rules~\cite{moin_2010}. 
Since the use of fixed quadrature points may produce overfitting during training~\cite{JANDER_integration}, we consider piecewise-linear solutions over a mesh.

The accuracy of mesh-based solutions depends on the mesh used to numerically solve the PDE \cite{adaptatividad_origienes,adaptatividad_origienes1,adaptatividad_origienes2}.
DNNs can also be employed to adapt the mesh.
In particular, \cite{Alfonzetti} proposes a mesh generator based on NNs to grow the number of initial elements to a user-selected one. 
The algorithm is based on the node probability density function obtained from an error estimate.
The authors of \cite{MANEVITZ2005447} propose the use of a DNN to obtain an effective mesh refinement in time-dependent PDEs solved by the finite-element method (FEM). They predict the areas of ``interest’’ at the next time stage from the information of previous times. 
In the same line, \cite{BOHN202161} employs recurrent DNNs to refine the mesh. These works use DNNs to adapt the mesh but they resort to traditional methods as FEM or finite differences to find the PDE solution.

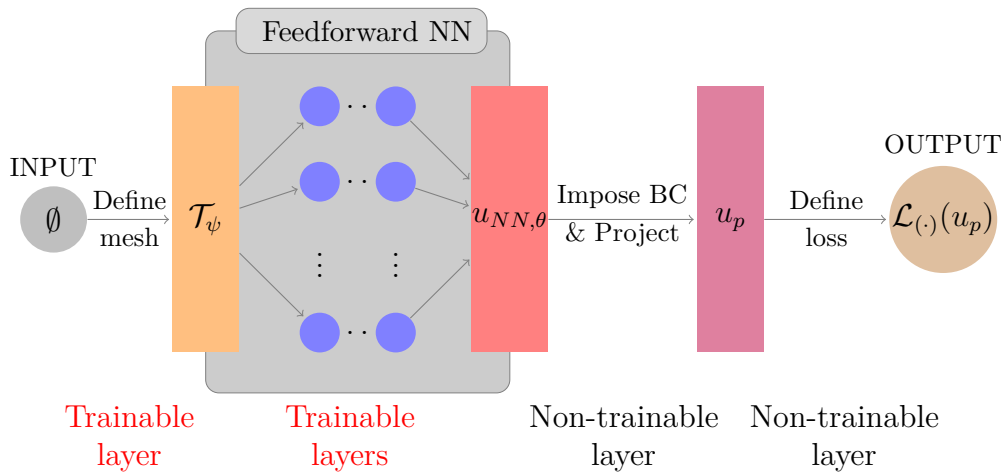
\begin{figure}[!ht]
  \centering
\def\layersep{1cm}
\begin{tikzpicture}[ shorten >=1pt,->,draw=black!50, node distance=\layersep]
    \tikzstyle{every pin edge}=[<-,shorten <=1pt]
    \tikzstyle{neuron}=[circle,fill=black!25,minimum size=25pt,inner sep=0pt]
    \tikzstyle{neuron2}=[circle,fill=black!25,minimum size=15pt,inner sep=0pt]
    \tikzstyle{mesh lay}=[rectangle,fill=orange!50,minimum width=25pt, minimum height=100pt,inner sep=0pt]
	\tikzstyle{hidden neuron}=[neuron2, fill=blue!50];
	\tikzstyle{ghost neuron}=[neuron, fill=gray!40!white, , opacity=1];
    \tikzstyle{empty neuron}=[neuron, fill=gray!50!white, , opacity=1];
    \tikzstyle{output lay}=[rectangle,fill=red!50,minimum width=25pt, minimum height=100pt,inner sep=0pt]
    \tikzstyle{multi lay}=[rectangle,fill=purple!50,minimum width=25pt, minimum height=100pt,inner sep=0pt]		
	\tikzstyle{final neuron}=[neuron, fill=brown!50];
\draw[very thin, rounded corners=5pt,fill=gray!40!white ,opacity=1] (0,1.5) rectangle (4,-3.3);
 \node[empty neuron] (empty) at (-2,-1) {$\emptyset$};
    \node[text width=4em, text centered, above of=empty, node distance=0.7cm] {\footnotesize INPUT};
	   
	\node[mesh lay, right of=empty,  xshift=1cm] (x) {$\mathcal{T}_{\psi}$};    

	\path[every node/.style={anchor=south}] (empty) edge node (Textempty) {\footnotesize Define} (x);
	\node [below of = Textempty, node distance=0.5cm] {\footnotesize mesh};
\path[yshift=1.5cm] node[ghost neuron, xshift=0.2cm] (Hghosta) at (0.7cm+\layersep,-1 cm) {$\hdots$};
\path[yshift=1.5cm] node[ghost neuron, xshift=0.2cm] (Hghostb) at (0.7cm+\layersep,-2 cm) {$\hdots$};
\path[yshift=1.5cm] node[ghost neuron, xshift=0.2cm] (Hghostc) at (0.7cm+\layersep,-4 cm) {$\hdots$};
\path[yshift=1.5cm] node[hidden neuron, xshift=0.5cm] (H1_1) at (\layersep,-1 cm) {};	
\path[yshift=1.5cm] node[hidden neuron, xshift=0.5cm] (H1_2) at (\layersep,-2 cm) {};
\path[yshift=1.5cm] node[ghost neuron, xshift=0.5cm] (Hghost1) at (\layersep,-3 cm) {$\vdots$};
\path[yshift=1.5cm] node[hidden neuron, xshift=0.5cm] (H1_3) at (\layersep,-4 cm) {};
\path[yshift=1.5cm] node[hidden neuron, xshift=0.0cm] (H2_1) at (1.5cm+\layersep,-1 cm) {};	
\path[yshift=1.5cm] node[hidden neuron, xshift=0.0cm] (H2_2) at (1.5cm+\layersep,-2 cm) {};
\path[yshift=1.5cm] node[ghost neuron, xshift=0.0cm] (Hghost2) at (1.5cm+\layersep,-3 cm) {$\vdots$};
\path[yshift=1.5cm] node[hidden neuron, xshift=0.0cm] (H2_3) at (1.5cm+\layersep,-4 cm) {};
    \node[output lay, right of=Hghost2, yshift=0.5cm, xshift=0.5cm] (O) {$u_{NN,\theta}$};   
	\node[multi lay, right of=O,  xshift=1.9cm] (M) {$u_{p}$};	
	\node[final neuron, right of=M, xshift=1.8cm] (F) {$\mathcal{L}_{(\cdot)}(u_{p})$};
    \foreach \dest in {1,...,3}
            \path (x) edge (H1_\dest);
    \foreach \source in {1,...,3}
        \path (H2_\source) edge (O);        
	\path[every node/.style={anchor=south}] (O) edge node (TextBC) {\footnotesize Impose BC} (M);
	\node [below of = TextBC, node distance=0.5cm] {\footnotesize \& Project};	
	\path[every node/.style={anchor=south}] (M) edge node (Textloss){\footnotesize Define} (F);
	\node [below of = Textloss, node distance=0.5cm] {\footnotesize loss};
\node[text width=4em, text centered, ,above of=F, node distance=1.cm] {\footnotesize OUTPUT}; 
\node[text width=6.5em, text centered] at (-1., -3.9) {\textcolor{red}{Trainable layer}}; 
\node[text width=6.5em, text centered] at (1.9, -3.9) {\textcolor{red}{Trainable layers}}; 
\node[text width=6.5em, text centered] at (5.5, -3.9) {Non-trainable layer}; 
\node[text width=6.5em, text centered] at (8.4, -3.9) {Non-trainable layer}; 
\draw[very thin, rounded corners=5pt,fill=gray!30!white ,opacity=1] (0.4,1.6+0.2) rectangle (3.7,1.+0.2);
\node [anchor = west]  at (0.6,1.3+0.2){\small Feedforward NN};    
    
\end{tikzpicture}
\caption{Architecture sketch of our proposed $r-$adaptive DNN.} \label{fig:sketch_r_adaptivity}
\end{figure}

Herein, we propose a deep $r-$adaptive method to simultaneously optimize the PDE solution and the location of the nodes in the corresponding mesh.
In our method, once we have selected the degree of the piecewise-polynomial approximation space and the number of elements in the mesh, we define a set of trainable parameters $\bs \psi$ that correspond to the location of the mesh nodes. 
From these variables, we define the mesh $\mathcal{T}_{\psi}$, which is used as an input of a DNN whose output is the piecewise-polynomial solution $u_p$. 
Then, the minimization of the loss function $\mathcal{L}$ automatically and simultaneously adapts the mesh $\mathcal{T}_{\psi}$ and the corresponding solution $u_p$.
A sketch of this architecture is shown in Figure \ref{fig:sketch_r_adaptivity}.

Our proposed method solves the optimal location of the nodes without solving auxiliary non-linear mesh moving PDE (MM-PDE) typical of traditional $r-$adaptive methods \cite{budd_huang_russell_2009,DORFI1987175,huang2011adaptive,MMPDE,Budd2}.
The technique of automatic-differentiation \cite{Autodiff,Autodiff_1,Autodiff_2} used within the back-propagation algorithm \cite{backprop} allows us to minimize directly the loss function.
Another substantial difference with traditional $r-$adaptive methods is that we can modify the topology of the mesh during the adaptation as long as we have access to the derivative of the algorithm that modifies the mesh with respect to each node location. 
In particular, when restricted to tensor product meshes, it is possible to fix some nodes in the mesh ---to ensure conformity with the material properties--- while allowing the others to freely move along the mesh.

This work shows the behavior of the $r-$adaptive method in combination with DRM, DLS, and PINNs to solve simple PDEs.
For the sake of simplicity, we use piecewise-linear functions over one-dimensional (1D) and two-dimensional (2D) tensor-product meshes to show the potential of the proposed method. 
The extension to higher-order piecewise-polynomials is straightforward.
We can also parametrize more complex geometries applying a mapping from a reference tensor-product mesh to the geometry of the physical space \cite{mesh}, possibly leaving some elements outside the computational domain (e.g., performing techniques similar to the Finite Cell Method \cite{cell_method,cell_method_2}). 
However, if we want to break the tensor product structure of the mesh and allow arbitrary changes in topology, it would require a mesh generator whose derivatives with respect to the node locations are accessible. 
The computation of the derivatives needed by general mesh-generation algorithms (like Delaunay triangularization \cite{PaulChew1989}) in TensorFlow~\cite{tensorflow2015-whitepaper} requires further future study. 
However, to overcome this restriction, projects like Autograd \cite{maclaurin2015autograd} (and its evolution JAX \cite{jax_github}) are progressing towards automatically differentiate native Python and NumPy functions.

The remainder of the manuscript is organized as follows. Section \ref{sc:ModelProblem} introduces the model problems and summarizes the numerical methods selected to solve them. Then, section \ref{sc:rMethod} introduces the $r-$adaptive method and explains how we incorporate it within a DNN.
Section \ref{sc:Implementation} is devoted to implementation details.
Section \ref{sc:Nresults} shows the numerical results. 
Finally, the last section summarizes the main conclusions and future work.

\section{Model problems and numerical methods}\label{sc:ModelProblem}

We consider two problems: a first-order hyperbolic problem and a second-order scalar elliptic problem. For the first (hyperbolic) problem, we introduce two losses, leading to a collocation and a least-squares method. For our second (elliptic) problem, we introduce the Deep Ritz method.

\subsection{Hyperbolic model problem}\label{sc:Hyperbolic}
Let $\Omega \in \mathbb{R}^d$ be a Lipschitz domain, where $d$ is the spatial dimension, and $\bs \beta(x)$ the velocity field.
We define the inflow and outflow boundary parts as ${\Gamma_- := \{ x \in \partial \Omega: \bs \beta \cdot \bs n < 0 \}}$ and  ${\Gamma_+ := \{x \in \partial \Omega: \bs \beta \cdot \bs n > 0 \}}$, respectively, with ${\partial \Omega = \Gamma_- \cup \Gamma_+}$, and $\bs n$ the outward normal vector to $\partial \Omega$.
We consider the boundary value problem, governed by the advection-reaction equation:
\begin{equation}\label{ec:advection}
\left\{
\begin{array}{rrcl}
 \nabla \cdot (\bs \beta \, u) + u(x) & = & f & \qquad \text{in} \quad \Omega,\\
u & = & 0 & \qquad \text{on} \quad \Gamma_-,
\end{array}
\right.
\end{equation}
where $f \in L^2(\Omega)$ is the source term and $u\in H_0^1(\Omega) = \{ u \in H^1(\Omega) \, | \, u=0 \text{ on }  \Gamma_- \}$ is the solution of the PDE.

Let $r := \nabla \cdot (\bs \beta \, u) + u - f$
be the residual. We define the loss for the collocation method as the $L^1-$norm of the residual, i.e.:
\begin{equation}\label{eq:losscol}
\mathcal{L}_{col} := ||r||_{L^1(\Omega)} = \int_{\Omega} |\nabla \cdot (\bs \beta \, u) + u - f|.
\end{equation}

We also introduce the loss associated to the $L^2$-norm of the residual (least-squares method):

\begin{equation}\label{eq:lossLS}
\mathcal{L}_{LS} := ||r||_{L^2(\Omega)} = \sqrt{\int_{\Omega} \left(\nabla \cdot (\bs \beta \, u) + u - f \right)^2}.
\end{equation}

The goal is to solve problem (\ref{ec:advection}) via minimization of the corresponding loss:
\begin{equation}\label{eq:minRitz}
u^* = \argmin_{u \in H_0^1(\Omega)} \mathcal{L}_{(\cdot)} (u).
\end{equation}

\subsection{Elliptic model problem} \label{sc:Elliptic}
We consider the open bounded domain $\Omega \in \mathbb{R}^{d}$, where $d$ is the spatial dimension, and we solve the following boundary value problem:
\begin{equation}
\left\{
\begin{array}{rrcl}
-\nabla \cdot (\sigma \nabla u) & = & f & \qquad \text{in} \quad \Omega, \\ 
u & = & 0 & \qquad \text{on} \quad \Gamma_D,\\  
(\sigma\nabla u) \cdot \bs{n}  & = & g & \qquad \text{on} \quad \Gamma_N,
\end{array}
\right.
\label{eq:no_poisson}
\end{equation}
where $\sigma \equiv \sigma(\bs x)\in L^ \infty(\Omega)$ are material properties,
${u\equiv u(\bsx) \in H_0^1(\Omega) = \{ u \in H^1(\Omega) \, | \, u=0 \text{ on }  \Gamma_D \}}$  is the solution of the second-order elliptic PDE, and $f \equiv f(\bsx) \in L^2(\Omega)$ is the volumetric source term. $\Gamma_D$ and $\Gamma_N$ denote the Dirichlet and Neumann parts of the boundary, respectively, with $\partial \Omega = \Gamma_D \cup \Gamma_N$,
$g \equiv g(\bsx)\big|_{\Gamma_N}\in H^{-1/2}(\Omega)$, and $\bs n$ denotes the unitary outward vector to $\partial \Omega$.
The energy function associated to the symmetric positive definite problem (\ref{eq:no_poisson}) is
\begin{equation}\label{eq:Ritz}
\mathcal{L}_{Ritz} = \frac{1}{2} \int_{\Omega} \sigma (\nabla u)^2 - \int_{\Omega} f u - \int_{\Gamma_N} g u.
\end{equation}
The Ritz method \cite{Ritz} finds the solution $u$ by minimizing $\mathcal{L}_{Ritz}$ following
Equation (\ref{eq:minRitz}).
As proven in \cite{Claes}, the solution $u$ of the boundary value problem (\ref{eq:no_poisson}) is also the solution of the minimization problem (\ref{eq:minRitz}).

\section{Deep $r-$adaptive method}\label{sc:rMethod}

We divide the proposed method into three parts, following the sketch shown in Figure \ref{fig:sketch_r_adaptivity}. The first one determines the location of the points that define the $r$-adaptive mesh. The second one uses a DNN to approximate the solution of the PDE at the mesh nodes, imposes the Dirichlet BCs, and interpolates the solution into the selected piecewise-polynomial space. The third part states the minimization problem and show how to compute the loss function. Finally, we explain our minimization algorithm.

\subsection{Adaptive mesh}
In this work, we discretize the $d$-dimensional domain $\Omega$ by a tensor product mesh $\mathcal{T_{\psi}}$. 
For each dimension, we build an ordered one-dimensional vector $\bs x$.
This vector contains fixed and variable coordinates. 
Fixed coordinates define the boundary and interfaces to separate materials. 
The variable coordinates are optimized during the training process to adapt the mesh.

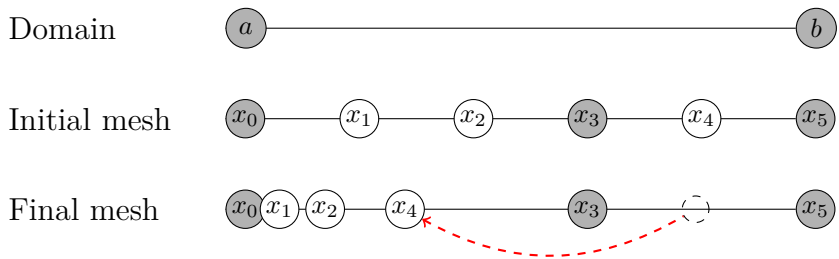
\begin{figure}[!htb]
  \centering
\begin{tikzpicture}
\node[anchor=north west] at (0,1.2) {
  \begin{tikzpicture}[x=1.5cm, y =1cm]  
  \node[anchor=west](a) at (-2,0){Domain};
	\def\nNodes{5}		
	\draw[-](0,0)--(\nNodes,0);
	\node [anchor = west,circle, inner sep=1pt, fill=gray!60!white, draw=black, minimum size=15pt] (x0) at (0,0) {\footnotesize $a$};
	\node [anchor = west,circle, inner sep=1pt, fill=gray!60!white, draw=black, minimum size=15pt] (x5) at (5,0) {\footnotesize $b$};	
  \end{tikzpicture}
  };

\node [anchor=north west] at (0,0) {
  \begin{tikzpicture}[x=1.5cm, y =1cm]  
  \node[anchor=west](a) at (-2,0){Initial mesh};
	\def\nNodes{5}		
	\draw[-](0,0)--(\nNodes,0);
	\node [anchor = west,circle, inner sep=1pt, fill=gray!60!white, draw=black, minimum size=10pt] (x0) at (0,0) {\footnotesize $x_0$};
	\node [anchor = west,circle, inner sep=1pt, fill=white, draw=black, minimum size=10pt] (x1) at (1,0) {\footnotesize $x_1$};
	\node [anchor = west, circle, inner sep=1pt, fill=white, draw=black, minimum size=10pt] (x2) at (2,0) {\footnotesize $x_2$};
	\node [anchor = west,circle, inner sep=1pt, fill=gray!60!white, draw=black, minimum size=10pt] (x3) at (3,0) {\footnotesize $x_3$};
	\node [anchor = west,circle, inner sep=1pt, fill=white, draw=black, minimum size=10pt] (x4) at (4,0) {\footnotesize $x_4$};
	\node [anchor = west,circle, inner sep=1pt, fill=gray!60!white, draw=black, minimum size=10pt] (x5) at (5,0) {\footnotesize $x_5$};	
  \end{tikzpicture}
  };

\node [anchor=north west] at (0,-1.2) {
 \begin{tikzpicture}[x=1.5cm, y =1cm]  
	\def\nNodes{5}	
  \node[anchor= west](a) at (-2.,0){Final mesh}; 
	\draw[-](0,0)--(\nNodes,0);
	\node [anchor = west,circle, inner sep=1pt, fill=gray!60!white, draw=black, minimum size=10pt] (x0) at (0,0) {\footnotesize $x_0$};
	\node [anchor = west,circle, inner sep=1pt, fill=white, draw=black, minimum size=10pt] (x1) at (0.3,0) {\footnotesize $x_1$};
	\node [anchor = west,circle, inner sep=1pt, fill=white, draw=black, minimum size=10pt] (x2) at (0.7,0) {\footnotesize $x_2$};
	\node [anchor = west,circle, inner sep=1pt, fill=white, draw=black, minimum size=10pt] (x4) at (1.4,0) {\footnotesize $x_4$};	
	\node [anchor = west,circle, inner sep=1pt, fill=gray!60!white, draw=black, minimum size=10pt] (x3) at (3,0) {\footnotesize $x_3$};
	\node [anchor = west,circle, inner sep=1pt, fill=gray!60!white, draw=black, minimum size=10pt] (x5) at (5,0) {\footnotesize $x_5$};	
	\node [anchor = west,circle, dashed, inner sep=1pt, draw=black, minimum size=10pt] (x4old) at (4,0) {};
	\draw [<-,draw = red, dashed, thick,] (x4) to [bend left=-30] (x4old);	
  \end{tikzpicture} 
};
\end{tikzpicture}
\caption{Representation of a one-dimensional domain, as well as the initial and final meshes. The mesh is built with one interior fixed node and three variables nodes. The movement of the node $x_4$ changes the topology in the final mesh, requiring to reorder the nodes.} \label{fig:topology}
\end{figure}

We allow changes in the mesh topology, as illustrated in Figure \ref{fig:topology}. 
In these, we observe that the trainable coordinates are allowed to \textit{jump} from one side to the other of the fixed node  $x_3$. 
The coincidence of coordinates at nodes is also allowed. 
If this scenario occurs, elements with zero volume are generated and have no contribution to the loss function.

In this training block, the set of trainable parameters $\bs \psi$ defines the location of the mesh nodes.
Then, we use a meshing algorithm ---trivial for tensor product meshes--- to generate the mesh $\mathcal{T_{\psi}}$.

\subsection{PDE solution approximation}\label{sc:PDE}

We approximate the solution $u$ of problems (\ref{eq:no_poisson}) and (\ref{ec:advection}) as
\begin{equation}
u \approx u_{p} := P \circ D \circ u_{NN_{\theta}}({\mathcal{T_{\psi}}}),
\end{equation}
where 
$u_{NN_{\theta}}$ is a DNN with a set of trainable parameters denoted as $\bs \theta$ that maps the mesh $\mathcal{T}_{\psi}$ into a surrogate solution $\tilde{u}$.
$D$ is a transformation used to impose Dirichlet BCs in a strong form.
Specifically,  $D(\tilde{u}):= u_D + \phi_D \tilde{u}$, where $\phi_D$ is a positive function in the interior of the domain and takes value zero on the Dirichlet boundary, and $u_D$ is a function to impose the lift for the case of inhomogeneous Dirichlet BC. 
Finally,
$P$ is the operator that builds the piecewise-linear solution over the mesh $\mathcal{T}_{\psi}$ by interpolating evaluations of the solution at mesh nodes.

\subsection{Minimization problem and loss function computation}

To find the approximate solution $u_{p}$ we rewrite the minimization problem from Eq. (\ref{eq:minRitz}):
\begin{equation}\label{eq:full_minimization}
u \approx u_{p}(\bs s^*) := \underset{
u_{p}(\bs s), \; \bs s \in S
} 
{\textup{ arg min }} 
\mathcal{L}_{(\cdot)}\left(u_{p}(\bs s)\right),
\end{equation}
where $\bs s$ denotes the set of all trainable parameters
\begin{equation}\label{ec_opt}
\bs s := \bigcup_{i, j=1}^{k_{\theta}, k_{\psi} } \left\lbrace  \theta_i, \psi_j  \right\rbrace,
\end{equation}
and $\theta_i$ are the trainable parameters used to approximate the PDE solution, while $\psi_j$ are the unknown node locations. Note that by construction, the search space $\{ u_{p}(\bs s), \bs s \in S\} \subset H_0^1$ is not a vector space in general.

We select a Gaussian quadrature rule to compute the integrals that appear in the loss functions ---equations (\ref{eq:losscol}), (\ref{eq:lossLS}), and (\ref{eq:Ritz})---. 
The choice of the number of quadrature points per element depends on the degree $p$ of the piecewise-polynomial function selected to approximate the solution and the source term function $f$ or the Neumman BC $g$ when it applies. 
Integration errors can negatively affect the optimization process, leading to poor approximation of the gradients used to update the DNN trainable parameters.

\subsection{Guided optimization}\label{sc:GuidedOPT}
 
The simultaneous minimization of the mesh and solution is a highly non-linear problem and may lead to inaccurate results due to the presence of local minima. 
For this reason, we decide to guide the optimization process in two stages.

The first stage only optimizes over the values theta, as in (\ref{ec_opt}), which corresponds to optimizing over the values of the DNN at the interpolation points over the fixed, initial mesh. 
In the second stage, we optimize both the solution and the mesh nodes location. This strategy not only increases the chances of avoiding local minima, but it also often prevents integration problems in the right-hand side $f$.
To distinguish both training phases, we will denote $u_p$ to the solution after optimization of the solution over a fixed grid. 
In contrast, $u_{p,r}$ will denote the solution after simultaneously optimizing both the solution and the node locations.


\section{Implementation}\label{sc:Implementation}
This section explains the algorithms used to define the tensor product mesh and the deep NN.

%
\subsection{Adaptive mesh implementation}

For each spatial dimension, we consider a vector of one-dimensional coordinates $\bs x$ that will be used to build a tensor product mesh.
We consider two types of coordinates: fixed and variable.
Among the fixed coordinates, we distinguish between those on the boundary, denoted as $a, b\in \partial \Omega$, with $a<b$, and the vector $\bs x^{fix}$ that stores $n_{fix}$ fixed coordinates between $a$ and $b$.
Incorporating interior fixed coordinates allows for the construction of conforming meshes needed to model problems with several materials.
On the other hand, the vector of variable coordinates $\bs x^{\delta}$ contains a set of $n_{\delta}$ coordinates, which are free to move during the training process, allowing for $r-$adaptation.

Algorithm \ref{alg:coord} computes a one-dimensional ordered vector $\bs x$ from $\bs x^{fix}$ and the vector of coordinates $\bs \psi$ updated during the training. 
$\bs \psi := \{a, \bs x^{\delta}, b\}$, and we need to map it into the computational domain after each iteration of the training.
The number of entries of $\bs x$ denoted
as $n := 2 + n_{fix} + n_{\delta}$ refers to the total number of mesh nodes along one dimension. 
We start the training from a uniform mesh initializing $\bs \psi$ as: $\psi_i \gets i, i=1$ to $n_{\delta}+2$.
The computational cost of the algorithm is $\mathcal{O}(n\log{}n)$. 

\begin{algorithm}
\small
 \KwInput{$\bs \psi, a, b, \bs x^{fix}, n_{fix}, n_{\delta}$} 
 \KwOutput{$\bs x$ (dim$(\bs x) =n = 2+n_{fix}+ n_{\delta}$)}

$\bs \psi \gets $sort($\bs \psi$) \tcp*{reorder $\bs \psi$}
 
$\bs{\psi} \gets (\bs \psi - \psi_1) / (\psi_{n_{\delta}+2} - \psi_1) (b - a) + a$ \tcp*{scale $\psi$ to the computational domain}
\uIf{$n_{fix}= 0$}{
    $\bs x \gets \bs{\psi} $\;
    }
\Else{
$\bs x \gets$ sort(concatenate($\bs{\psi}, \bs x^{fix}$)) \tcp*{add interior fixed coordinates and reorder}
  }

 \Return{$\bs x$}\;
 \caption{1D computation of the coordinate vector: Training coordinates}
 \label{alg:coord}
\end{algorithm}

Finally, we construct the tensor product mesh $\mathcal{T_{\psi}} $ by the combination of the $d$ one-dimensional vectors $\bs x$. 
The number of nodes $n_{\mathcal{T}}$ and elements $n_K$ of  $\mathcal{T_{\psi}}$ is
\begin{equation}
n_{\mathcal{T}}:= \prod_{i=1}^d n_i, \quad n_K:= \prod_{i=1}^d (n_i -1),
\end{equation}
where $n_i$ denotes the dimension of the vector $\bs x$ associated to the spatial dimension $i$.


\subsection{Deep Neural Network}\label{sc:DNN}
Once we have the $n_{\mathcal{T}}$ points in the physical space that define the mesh $\mathcal{T_{\psi}}$, we use a feed-forward NN $u_{NN\theta}(\mathcal{T_{\psi}}): \mathbb{R}^{d} \rightarrow \mathbb{R}$ composed of $k$ layers as follows:
\begin{equation}\label{ec:feedforward}
u_{NN,\theta}:= \ell^{(k)} \circ \ell^{(k-1)} \circ \dots \circ \ell^{(1)}(\mathcal{T_{\psi}}).
\end{equation}
Each layer $\ell^{(i)}$ is a non-linear mapping composed of an activation function $\alpha^{(i)}$ applied component-wise to an affine transformation,
\begin{equation}
\ell^{(i)}(\cdot) := \alpha^{(i)} \left( \bs W^{(i) \top} (\cdot) + \bs b^{(i)} \right).
\end{equation}
In this work, we take our activation function to be the sigmoid function $\alpha^{(i)}=\frac{1}{1+e^{-x}}$ if $i\neq k$, and the identity if $i=k$.
$N_i$ is the dimension of the layer output $\ell^{(i)}$ and is usually understood as the number of neurons in the layer. It is related to the number of trainable parameters, the so-called weights $\bs W^{(i)} \in \mathbb{R}^{N_i \times N_{i-1}}$ and biases $\bs b^{(i)}\in \mathbb{R}^{N_i}$ of the affine transformation. 
We denote the set of all trainable parameters of the DNN $u_{NN,\theta}$ as:
\begin{equation}
\displaystyle \bs \theta:= \bigcup_{i=1}^{k} \left\lbrace \left( \bs W^{(i)}, \bs b^{(i)} \right) \right\rbrace.
\end{equation}

\section{Numerical results}\label{sc:Nresults}

We adopt the strategy described in Section \ref{sc:GuidedOPT} for training. We first optimize the solution over a fixed mesh (denoted as $u_p$), then we optimize both the mesh and solution to obtain $u_{p,r}$.
In both stages, we use ADAM optimizer \cite{adam}.
The selected feed-forward NN has five hidden layers of ten neurons, each activated with a sigmoid function ---see Section \ref{sc:DNN} for details.

Experiment 1 solves the hyperbolic problem with the $r-$adaptive method in combination with residual methods, as explained in Section \ref{sc:Hyperbolic}.
Experiments 2 to 6, summarized in Table \ref{tab:experiments}, focus on solving the elliptic problem introduced in Section \ref{sc:Elliptic} by the Deep Ritz Method of Eq.~\ref{eq:Ritz}.

\begin{table}[htb!]
   \centering
   \caption{Numerical experiments for the elliptic problem of Section \ref{sc:Elliptic}: definition and exact solutions. }
   \label{tab:experiments}
   \begin{tabular}{cccc}
     \hline
\begin{tabular}{c}Experiment \\  ID\end{tabular}         &  $u_{exact}$   & $\mathcal{L}_{Ritz}(u_{exact}) $ & \begin{tabular}{c} Domain, material \\ properties, and BCs \end{tabular} \rule{0pt}{5ex}\rule[-3ex]{0pt}{0pt} \\  
     \hline
 Ex. 2 &  $x^{0.7}$  &  -1.5385 & 
\begin{tabular}{c}
\begin{tikzpicture}
\node (A) at (-0.5, 0){};
\node (C) at (1.5, 0){};
\node (B) at (3.5, 0){};
\draw[thick,black] (A.center) to (B.center);
\node [circle, fill=blue, draw=blue, inner sep=0pt, minimum size=4pt] at (A) {};
\node [circle, fill=red, draw=red, inner sep=0pt, minimum size=4pt] at (B) {};
\node[anchor = west, red] at (B){$\Gamma_N$};
\node[anchor = east, blue] at (A){$\Gamma_D$};
\node[text centered, above of = C, node distance=0.2cm,xshift=0cm] {$\sigma = 1$};
\node[text centered, above of = C, node distance=-0.4cm,xshift=0cm] {$x \in (0,10)$};
\end{tikzpicture}
\end{tabular}
\rule{0pt}{5ex}\rule[-2ex]{0pt}{0pt} 
\\ 
Ex. 3 &  $\text{atan}(2x-10)+\text{atan}(10)$  &  -1.5701 &
\begin{tabular}{c}
\begin{tikzpicture}
\node (A) at (-0.5, 0){};
\node (C) at (1.5, 0){};
\node (B) at (3.5, 0){};
\draw[thick,black] (A.center) to (B.center);
\node [circle, fill=blue, draw=blue, inner sep=0pt, minimum size=4pt] at (A) {};
\node [circle, fill=red, draw=red, inner sep=0pt, minimum size=4pt] at (B) {};
\node[anchor = west, red] at (B){$\Gamma_N$};
\node[anchor = east, blue] at (A){$\Gamma_D$};
\node[text centered, above of = C, node distance=0.2cm,xshift=0cm] {$\sigma = 1$};
\node[text centered, above of = C, node distance=-0.4cm,xshift=0cm] {$x \in (0,10)$};
\end{tikzpicture}
\end{tabular}
\rule{0pt}{5ex}\rule[-2ex]{0pt}{0pt} 
\\ 
Ex. 4 &  $\dfrac{\sin (2 \pi x)}{\sigma}$  &  -5.8724 & 
\begin{tabular}{c}
\begin{tikzpicture}
\node (A) at (-0.5, 1.5){};
\node (C) at (1.5, 1.5){};
\node (B) at (3.5, 1.5){};

\draw[thick,black] (A.center) to (C.center);
\draw[thick,violet] (C.center) to (B.center);
\node [circle, fill=blue, draw=blue, inner sep=0pt, minimum size=4pt] at (A) {};
\node [circle, fill=blue, draw=blue, inner sep=0pt, minimum size=4pt] at (B) {};
\node [circle, fill=black, draw=black, inner sep=0pt, minimum size=2pt] at (C) {};

\node[anchor = west, blue] at (B){$\Gamma_D$};
\node[anchor = east, blue] at (A){$\Gamma_D$};

\node[text centered, above of = A, node distance=-0.4cm,xshift=0.cm] {$x =0$};
\node[text centered, above of = C, node distance=-0.4cm,xshift=0.0cm,] {$x =0.5$};
\node[text centered, above of = B, node distance=-0.4cm,xshift=0.0cm] {$x =1$};
\node[text centered, above of = A, node distance=0.2cm,xshift=1.0cm] {$\sigma=1$};
\node[text centered, above of = C, node distance=0.2cm,xshift=1.0cm, violet] {$\sigma=10$};
\end{tikzpicture}
\end{tabular}
\rule{0pt}{5ex}\rule[-2ex]{0pt}{0pt} 
\\
Ex. 5 &  $x_1^2 (x_1 -1)\, x_2^2 (x_2 -1)$  &  -0.0013 & 
\begin{tabular}{c}
\begin{tikzpicture}[x=1.25cm, y =1.25cm]
\draw[->](0,0)--(3,0);
\draw[->](0,0)--(0,3);
\node [anchor = south] at (3,0){$x_1$};
\node [anchor = east] at (0,3){$x_2$};
\node [anchor = east, blue] at (0,1.25){$\Gamma_D$};
\node [anchor = north, blue] at (1.25, 0){$\Gamma_D$};
\node [anchor = west, blue] at (2.5, 1.25){$\Gamma_D$};
\node [anchor = south, blue] at (1.25, 2.5){$\Gamma_D$};
\node [anchor = north ,] at (0,0){$(0,0)$};
\node [anchor = north  ,] at (2.5, 0){$(1,0)$};
\node [anchor =  east,] at (0, 2.5){$(0,1)$};
\node [anchor = west ,] at (2.5, 2.5){$(1,1)$};
\draw[very thick, blue, fill=gray!20!white ] (0,0) rectangle (2.5,2.5);
\node [anchor = north]at (1.25,1.5){$\sigma= 1 $};
\end{tikzpicture}
\end{tabular}
\rule{0pt}{5ex}\rule[-2ex]{0pt}{0pt} 
\\
Ex. 6 &  $r^{\frac{2}{3}} \sin \left(\frac{2}{3}  \left( \theta + \frac{\pi}{2} \right) \right)$  &  -0.9181 & 
\begin{tabular}{c}
\begin{tikzpicture}[x=1.25cm, y =1.25cm]
\draw[->](0,0)--(3,0);
\draw[->](0,0)--(0,3);
\node [anchor = south] at (3,0){$x_1$};
\node [anchor = east] at (0,3){$x_2$};
\path[draw=black,fill=gray!20!white] (1.25,1.25)--(0,1.25)--(0,2.5)--(2.5,2.5)--(2.5,0)--(1.25,0)--cycle;
\node [anchor = north east, blue] at (1.25,1.25){$\Gamma_D$};
\node [anchor = west, red] at (2.5, 1.25){$\Gamma_N$};
\node [anchor = south, red] at (1.25, 2.5){$\Gamma_N$};
\node [anchor = east, red] at (0, 1.875){$\Gamma_N$};
\node [anchor = south, red] at (1.875,0){$\Gamma_N$};
\node [anchor =  east,] at (0,1.25){$(-1,0)$};
\node [anchor = north  ,] at (2.5, 0){$(1,-1)$};
\node [anchor = north  ,] at (1.25, 0){$(0,-1)$};
\node [anchor =  east,] at (0, 2.5){$(-1,1)$};
\node [anchor = west ,] at (2.5, 2.5){$(1,1)$};
\node [anchor = north]at (0.75,2.25){$\sigma= 1 $};
\draw [fill=green!30] (1.25,1.25) -- ++(1cm,0) arc(0:45:1cm) node[midway,left] {$\theta$} -- cycle;
\draw[thick, ->]  (1.25,1.25) -- node[above] {$r$} (2.1,2.1);
\draw[thin,dashed ]  (2.1,1.25) -- (2.5,1.25);
\path[draw=blue,very thick,] (0,1.25)--(1.25,1.25)--(1.25,0);
\path[draw=red,very thick,](0,1.25)--(0,2.5)--(2.5,2.5)--(2.5,0)--(1.25,0);
\end{tikzpicture}
\end{tabular}
\rule{0pt}{5ex}\rule[-2ex]{0pt}{0pt}
\\
\end{tabular}
\end{table}

\subsection{Experiment 1: Residual methods}

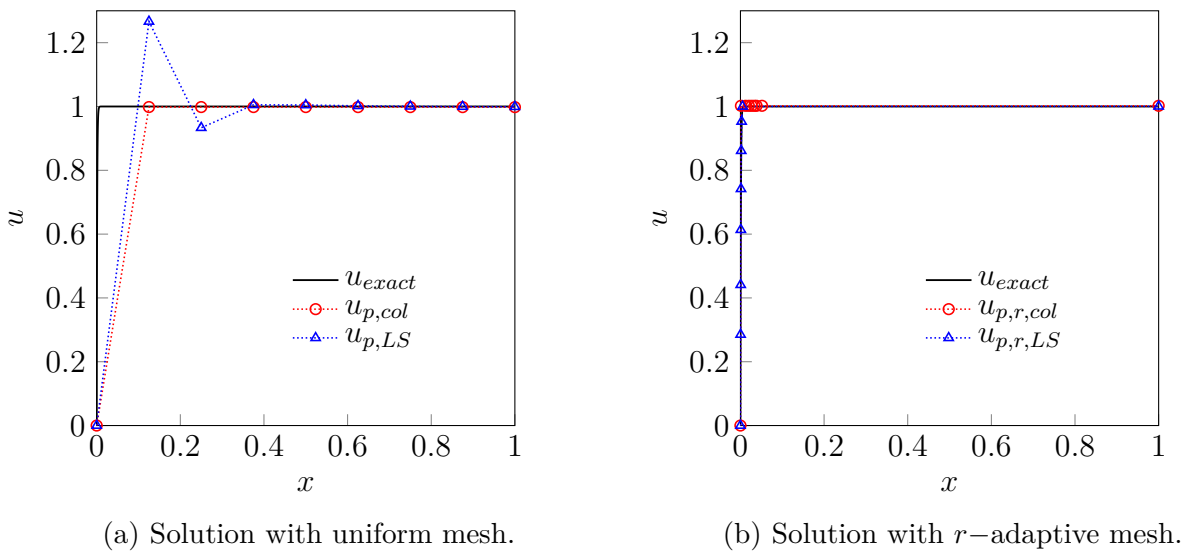
\begin{figure}[htb!]
  \centering
  \begin{subfigure}{0.49\textwidth}
	\begin{tikzpicture}
	\begin{axis}[scale only axis, xlabel = $x$, ylabel = $u$, 
	ytick pos=left,
	xtick pos=left,
	y label style={at={(-0.15,0.5)}}, 
	height=5.5cm, width=5.5cm, 
	xmin=0, ymin=0, 
	xmax=1, ymax=1.3,
	legend columns = 1,
	legend style= {at={(0.8,0.4)},draw=none,fill=none,nodes={scale=1, transform shape}}, 
	legend cell align={left}
	]
	%
	\addplot [black, line width=0.7pt,domain=0:10] table[x=x ,y=u_exact]{data/beta_0.001/r_L2/exact.csv};
	\addlegendentry{$u_{exact}$ };
	\addplot[line width=0.6pt,color=red, densely dotted,mark=o, mark options={solid}] table[x=x ,y=u_pred]{data/beta_0.001/r_L1/partition_0.csv};
	\addlegendentry{$u_{p,col}$ };
	\addplot[line width=0.6pt,color=blue, densely dotted,mark=triangle, mark options={solid}] table[x=x ,y=u_pred]{data/beta_0.001/r_L2/partition_0.csv};
	\addlegendentry{$u_{p,LS}$ };
	\end{axis}
\end{tikzpicture}
   \caption{Solution with uniform mesh. \label{fig:Ex_beta2_panel1}}
  \end{subfigure}
  \begin{subfigure}{0.49\textwidth}
	\begin{tikzpicture}
	\begin{axis}[scale only axis, xlabel = $x$, ylabel = $u$, 
	ytick pos=left,
	xtick pos=left,
	y label style={at={(-0.15,0.5)}}, 
	height=5.5cm, width=5.5cm, 
	xmin=0, ymin=0, 
	xmax=1, ymax=1.3,
	legend columns = 1,
	legend style= {at={(0.8,0.4)},draw=none,fill=none,nodes={scale=1, transform shape}}, 
	legend cell align={left}
	]
	%
	\addplot [black, line width=0.7pt,domain=0:10] table[x=x ,y=u_exact]{data/beta_0.001/r_L2/exact.csv};
	\addlegendentry{$u_{exact}$ };
	\addplot[line width=0.6pt,color=red, densely dotted,mark=o, mark options={solid}] table[x=x ,y=u_pred]{data/beta_0.001/r_L1/partition_1.csv};
	\addlegendentry{$u_{p,r,col}$ };
	\addplot[line width=0.6pt,color=blue, densely dotted,mark=triangle, mark options={solid}] table[x=x ,y=u_pred]{data/beta_0.001/r_L2/partition_1.csv};
	\addlegendentry{$u_{p,r,LS}$ };
	\end{axis}
	\end{tikzpicture}
   \caption{Solution with $r-$adaptive mesh. \label{fig:Ex_beta2_panel2}}
  \end{subfigure}
\caption{Experiment 1. Exact and approximate solutions for $\beta = 10^{-3}$ using uniform (left panel) and $r-$adaptive (right panel) meshes.} 
\label{fig:Ex_residual2}
\end{figure}

We solve the problem of Equation \ref{ec:advection} with $f=1$. The exact solution   is $u = 1-e^{\frac{-x}{\beta}}$.
Figure 
\ref{fig:Ex_residual2} shows the exact and approximated solutions for $\beta=10^{-3}$.
We consider a mesh of eight elements. The approximated solutions are computed via the minimization of the residuals $\mathcal{L}_{col}$ defined in Eq. (\ref{eq:losscol}) and $\mathcal{L}_{LS}$ in Eq. (\ref{eq:lossLS}). 
Figure \ref{fig:Ex_beta2_panel1} shows the Gibbs effect for the Least-Square solution $u_{p,LS}$. 
The Gibbs phenomenon is not present with the collocation method.
When relocating the nodes (see Figure \ref{fig:Ex_beta2_panel2}), we obtain superior-quality solutions that are free of Gibbs fluctuations.

\subsection{Experiment 2: Singular solution}
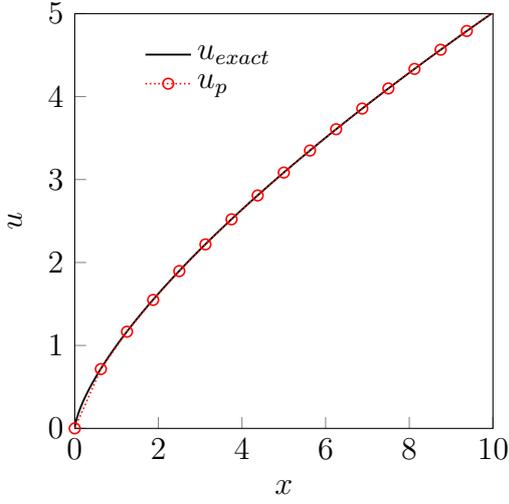
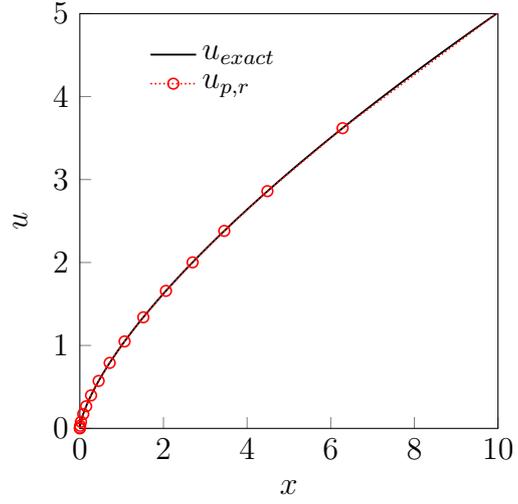
\begin{figure}[htb!]
  \centering
  \begin{subfigure}{0.49\textwidth}
	\begin{tikzpicture}
	\begin{axis}[scale only axis, xlabel = $x$, ylabel = $u$, 
	ytick pos=left,
	xtick pos=left,
	y label style={at={(-0.1,0.5)}}, 
	height=5.5cm, width=5.5cm, 
	xmin=0, ymin=0, 
	xmax=10, ymax=5,
	legend columns = 1,
	legend style= {at={(0.5,0.95)},draw=none,fill=none,nodes={scale=1, transform shape}}, 
	legend cell align={left}
	]
	%
	\addplot [black, line width=0.7pt,domain=0:10] table[x=x ,y=u_exact]{data/pb_121/exact.csv};
	\addlegendentry{$u_{exact}$ };
	\addplot[line width=0.6pt,color=red, densely dotted,mark=o, mark options={solid}] table[x=x ,y=u_pred]{data/pb_121/partition_0.csv};
	\addlegendentry{$u_{p}$ };
	\end{axis}
\end{tikzpicture}
   \caption{Solution with uniform mesh. \label{fig:test1_mesh}}
  \end{subfigure}
  \begin{subfigure}{0.49\textwidth}
	\begin{tikzpicture}
	\begin{axis}[scale only axis, xlabel = $x$, ylabel = $u$, 
	ytick pos=left,
	xtick pos=left,
	y label style={at={(-0.1,0.5)}}, 
	height=5.5cm, width=5.5cm, 
	xmin=0, ymin=0, 
	xmax=10, ymax=5,
	legend columns = 1,
	legend style= {at={(0.5,0.95)},draw=none,fill=none,nodes={scale=1, transform shape}}, 
	legend cell align={left}
	]
	%
	\addplot [black, line width=0.7pt,domain=0:10] table[x=x ,y=u_exact]{data/pb_121/exact.csv};
	\addlegendentry{$u_{exact}$ };
	\addplot[line width=0.6pt,color=red, densely dotted,mark=o, mark options={solid}] table[x=x ,y=u_pred]{data/pb_121/partition_1.csv};
	\addlegendentry{$u_{p,r}$ };
	\end{axis}
	\end{tikzpicture}
   \caption{Solution with $r-$adaptive mesh. \label{fig:test1_mesh_r}}
  \end{subfigure}
\caption{Experiment 2. Solutions exact and approximate at the end of each training step.} 
\label{fig:Test_1}
\end{figure}

The solution of Ex. 2  has a singularity at $x=0$, where the derivative is infinite. 
Figure \ref{fig:test1_mesh} shows the solution for a fixed 16-elements mesh, and Figure \ref{fig:test1_mesh_r} shows the final solution with our proposed $r-$adaptive method. 
We observe small elements accumulating nearby the singularity, creating a geometrical grid as expected for this type of singularity \cite{Strouboulis}.

\begin{figure}[htb!]
  \centering
  \begin{subfigure}{0.49\textwidth}
	\begin{tikzpicture}  
	\begin{axis}[scale only axis, xlabel = Epoch, ylabel = $\mathcal{L}_{Ritz}(u_{p,r}) - \mathcal{L}_{Ritz}(u_{exact})$, 
	ymode=log,
	ytick pos=left,
	xtick pos=left,
	y label style={at={(-0.17,0.5)}}, 
	height=5.5cm, width=5.5cm, 
	xmin=0, ymin=0.01, 
	xmax=8000, ymax=1,
	legend cell align={left}
	]	
	\addplot[color=black, mark options={solid}] table[x=epoch ,y=lossN]{data/pb_121/loss_clean1.csv};	
	\addplot[color=red, mark options={solid}] table[x=epoch ,y=lossN]{data/pb_121/loss_clean0.csv};	
	\draw[dashed, color=red](1000,0.001)--(1000,1);
	\node [rotate=90, color = red] at (500,0.3){\footnotesize 1st training step};	
	\node at (4000,0.7){\footnotesize 2nd training step};	
	\end{axis}
	\end{tikzpicture}  
   \caption{Loss evolution for 16 elements. \label{fig:test1_loss}}
  \end{subfigure}
  \begin{subfigure}{0.49\textwidth}
	\begin{tikzpicture}  
	\begin{axis}[scale only axis, xlabel = $\#$ of elements, ylabel = $||u-u_{(\cdot)}||_{H_{1(\Omega)}}$, 
	xmode=log,
	ymode=log,
	ytick pos=left,
	xtick pos=left,
	y label style={at={(-0.01,0.5)}}, 
	height=5.5cm, width=5.5cm, 
	xmin=1, ymin=0.1, 
	xmax=1024, ymax=1,
	log basis x={2},
	legend cell align={left}
	]
	\addplot[dashed,color=red, ,mark=o, mark options={solid}] table[x=elements ,y=static]{data/pb_121/norma.csv};
	\addplot[dashed,color=black, ,mark=diamond, mark options={solid}] table[x=elements ,y=static_FEM]{data/pb_121/norma.csv};
	\addplot[restrict x to domain=1:6, dotted,color=blue, ,mark=triangle, mark options={solid}] table[x=elements ,y=r]{data/pb_121/norma.csv};
	\addplot[restrict x to domain=1:6, dotted,color=black, ,mark=square, mark options={solid}] table[x=elements ,y=r_FEM]{data/pb_121/norma.csv};
	\node at (13,0.54){\footnotesize \textcolor{red}{${(p)}$}};
	\node at (13.5,0.35){\footnotesize \textcolor{black}{${(\text{FEM})}$}};
	\node at (40,0.2){\footnotesize \textcolor{blue}{${(p,r)}$ }};
	\node at (40,0.12){\footnotesize \textcolor{black}{${(\text{FEM},r)}$}};
	\addplot [black] coordinates {(128 ,0.35) (512,0.265238)(512,0.35)(128 ,0.35)} ;
	\node at (275,0.38){\footnotesize 1:5};
	\end{axis}
	\end{tikzpicture}  
   \caption{Convergence of the solutions in energy norm. \label{fig:test1_conv}}
  \end{subfigure}
\caption{Experiment 2. Loss evolution and convergence of the solutions} 
\label{fig:Test_1add}
\end{figure}
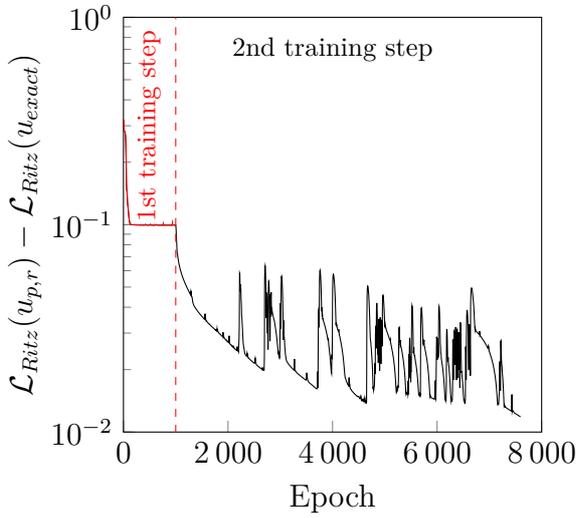
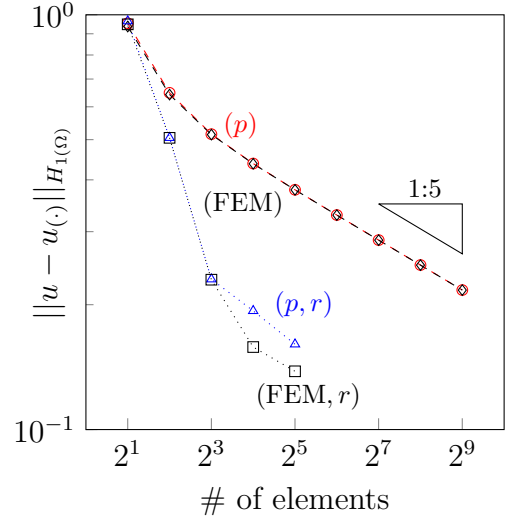

Figure \ref{fig:test1_loss} shows the error in the loss for a solution over a mesh of 16 elements. 
We observe a decrease in the loss error values during the second training step, i.e., when we optimize the node locations. 
This indicates the improvement in accuracy of the solution $u_{p,r}$ over $u_{p}$.

For different number of elements in the mesh, Figure \ref{fig:test1_conv} shows the energy-norm error, defined as: 
\begin{equation}
||u-u_{(\cdot)}||_{H^1(\Omega)} = \sqrt{\int_{\Omega} \left( (u_{exact} - u_{(\cdot)})^2 +  (\nabla u_{exact}- \nabla u_{(\cdot)})^2 \right)}.
\end{equation}
We denote as $u_{\text{FEM}}$ and $u_{\text{FEM},r}$ the FEM solutions over the uniform and $r-$adaptive meshes, respectively.
The convergence rate for uniform refinements is 0.2 \cite{Gui1986b}. Figure \ref{fig:test1_conv} confirms this behavior both for the FEM and DL solutions for uniform meshes. When considering $r-$adaptive meshes, the convergence rate dramatically increases (see Figure \ref{fig:test1_conv}). However, the convergence rate deteriorates at some point due to either the optimizer and/or considered architecture (see the difference between the FEM and DL solutions on the $r-$adaptive meshes for 16 and 32 elements).

\subsection{Experiment 3: High-gradient solution}\label{experimentS}
\begin{figure}[!htb]
  \centering
  \begin{subfigure}{0.49\textwidth}
	\begin{tikzpicture}
	\begin{axis}[scale only axis, xlabel = $x$, ylabel = $u$, 
	ytick pos=left,
	xtick pos=left,
	y label style={at={(-0.1,0.5)}}, 
	height=5.5cm, width=5.5cm, 
	xmin=0, ymin=0, 
	xmax=10, ymax=3,
	legend columns = 1,
	legend style= {at={(0.45,0.95)},draw=none,fill=none,nodes={scale=1, transform shape}}, 
	legend cell align={left}
	]
	%
	\addplot [black, line width=0.7pt,domain=0:10] table[x=x ,y=u_exact]{data/pb_131/exact.csv};
	\addlegendentry{$u_{exact}$ };
	\addplot[line width=0.6pt,color=red, densely dotted,mark=o, mark options={solid}] table[x=x ,y=u_pred]{data/pb_131/partition_1.csv};
	\addlegendentry{$u_{p,r}$ };
	\end{axis}
\end{tikzpicture}
   \caption{Solution on $r-$adaptive mesh. \label{fig:test2_sol}}
  \end{subfigure}
  \begin{subfigure}{0.49\textwidth}
	\begin{tikzpicture}  
	\begin{axis}[scale only axis, xlabel = Epoch, ylabel = $\mathcal{L}_{Ritz}(u_{p,r}) - \mathcal{L}_{Ritz}(u_{exact})$, 
	ymode=log,
	ytick pos=left,
	xtick pos=left,
	y label style={at={(-0.17,0.5)}}, 
	height=5.5cm, width=5.5cm, 
	xmin=0, ymin=0.001, 
	xmax=6000, ymax=1,
	legend cell align={left}
	]	
	\addplot[color=black, mark options={solid}] table[x=epoch ,y=lossN]{data/pb_131/loss_clean1.csv};	
	\addplot[color=red, mark options={solid}] table[x=epoch ,y=lossN]{data/pb_131/loss_clean0.csv};	
	\draw[dashed, color=red](1000,0.001)--(1000,1);
	\node [rotate=90, color = red] at (500,0.01){\footnotesize 1st training step};	
	\node at (3000,0.003){\footnotesize 2nd training step};	
	\end{axis}
	\end{tikzpicture}
   \caption{Loss evolution. \label{fig:test2_loss}}
  \end{subfigure}
\caption{Experiment 3. Exact and approximate solutions and loss evolution.} 
\label{fig:Test_2}
\end{figure}
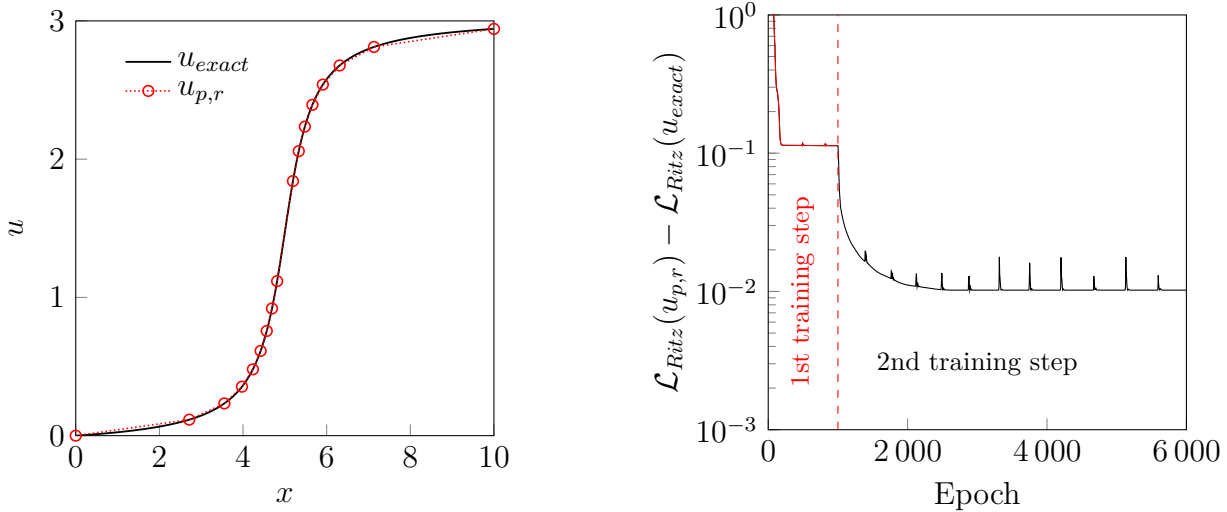

We now consider a solution with a strong gradient.
Figure \ref{fig:test2_sol} shows the exact and approximated solutions in a 16-element mesh. 
We observe an accumulation of elements in the zones with high gradients, as expected.
Figure \ref{fig:test2_loss} shows the loss evolution during the two steps of the training.
The decrease in the loss indicates that the solution obtained with the $r-$adaptive mesh is significantly better than the one calculated on a fixed uniform mesh.
  
In this experiment, the method diverges when executing only the one-step optimization without initializing the feed-forward NN on the uniform fixed mesh: 
The points accumulate in zones with small gradients. The quantity $\int_{\Omega} fu,$ increases, producing values of the loss $\mathcal{L}_{Ritz}$ lower than the best physically possible ---the one corresponding to the exact solution---, which indicates an integration error in the right-hand side. 
Although introducing some fixed points in the mesh solves the problem,  we decided to use the two-step training (see Section \ref{sc:GuidedOPT}) that also resolves the issue.

\subsection{Experiment 4: Discontinuous materials}

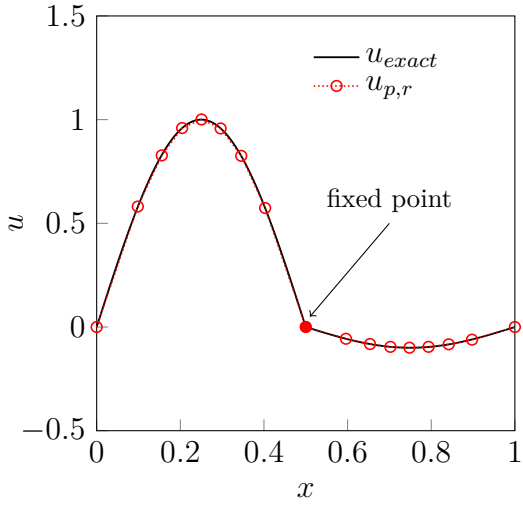
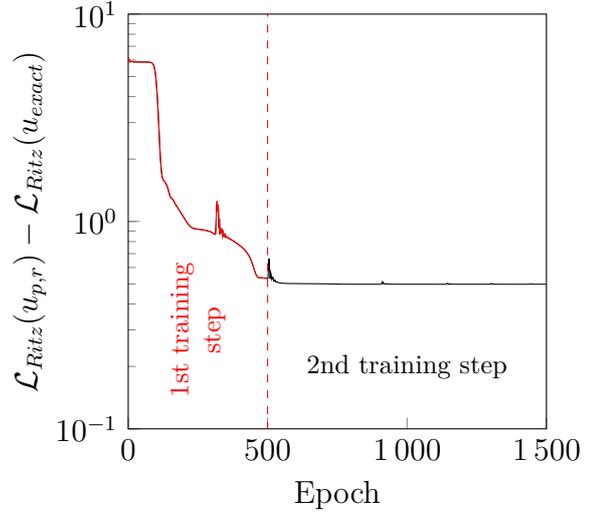
\begin{figure}[!htb]
  \centering
  \begin{subfigure}{0.49\textwidth}
	\begin{tikzpicture}
	\begin{axis}[scale only axis, xlabel = $x$, ylabel = $u$, 
	ytick pos=left,
	xtick pos=left,
	y label style={at={(-0.15,0.5)}}, 
	height=5.5cm, width=5.5cm, 
	xmin=0, ymin=-0.5, 
	xmax=1, ymax=1.5,
	legend columns = 1,
	legend style= {at={(0.85,0.95)},draw=none,fill=none,nodes={scale=1, transform shape}}, 
	legend cell align={left}
	]
	%
	\addplot [black, line width=0.7pt,domain=0:10] table[x=x ,y=u_exact]{data/pb_140/exact.csv};
	\addlegendentry{$u_{exact}$ };
	\addplot[line width=0.6pt,color=red, densely dotted,mark=o, mark options={solid}] table[x=x ,y=u_pred]{data/pb_140/partition_1.csv};
	\addlegendentry{$u_{p,r}$ };
	\addplot [only marks,mark=*,red,] coordinates {(0.5,0)};		
	\draw [<-,thin] (axis cs:0.51,0.05) -- (axis cs:0.7,0.5) node[above]{\footnotesize fixed point} ;
	
	\end{axis}
\end{tikzpicture}
   \caption{Solution on $r-$adaptive mesh. \label{fig:test3_sol}}
  \end{subfigure}
  \begin{subfigure}{0.49\textwidth}
	\begin{tikzpicture}  
	\begin{axis}[scale only axis, xlabel = Epoch, ylabel = $\mathcal{L}_{Ritz}(u_{p,r}) - \mathcal{L}_{Ritz}(u_{exact})$, 
	ymode=log,
	ytick pos=left,
	xtick pos=left,
	y label style={at={(-0.17,0.5)}}, 
	height=5.5cm, width=5.5cm, 
	xmin=0, ymin=0.1, 
	xmax=1500, ymax=10,
	legend cell align={left}
	]	
	\addplot[color=black, mark options={solid}] table[x=epoch ,y=lossN]{data/pb_140/loss_clean1.csv};	
	\addplot[color=red, mark options={solid}] table[x=epoch ,y=lossN]{data/pb_140/loss_clean0.csv};	
	\draw[dashed, color=red](500,0.001)--(500,10);
	\node [rotate=90, color = red] at (250,0.3){\footnotesize \begin{tabular}{c} 1st training \\ step\end{tabular} };	
	\node at (1000,0.2){\footnotesize 2nd training step};
	\end{axis}
	\end{tikzpicture}
   \caption{Loss evolution. \label{fig:test3_loss}}
  \end{subfigure}
\caption{Experiment 4. Exact and approximate solutions and loss evolution.} 
\label{fig:Test_3}
\end{figure}

This experiment reveals the capability of the proposed method to solve problems involving different materials. 
We include a fixed point to separate subdomains with different materials using a conforming mesh. 
Figure \ref{fig:Test_3} shows the solution and the evolution of the loss during the training process. 
Since the solution is smooth away from the fixed point, the improvement in the approximated solution accuracy between uniform mesh and the $r-$adaptive solution is lower than that observed in previous experiments.

\subsection{Experiment 5: two-dimensional (2D) smooth solution}

\begin{figure}[!htb]
  \centering
  \begin{subfigure}{0.49\textwidth}
  \begin{tikzpicture}[scale=0.9]
    \begin{axis}[view={0}{90},
    scaled z ticks=false,
	zticklabel=\pgfkeys{/pgf/number format/.cd,fixed,precision=2,zerofill}\pgfmathprintnumber{\tick}, 
    colorbar,
    xlabel = $x_1$,
    ylabel = $x_2$,
    height=7.cm, width=7cm, 
	xmin=0, ymin=0, 
	xmax=1, ymax=1,
	zmin = 0, zmax = 0.035,
	ztick = {0,0.01,0.02},  
    ]
      \addplot3[surf, point meta=explicit,] table [z expr=0.0, meta index=2] {data/pb_220/u_approx_1.csv};
      \addplot3[surf, point meta=explicit, mesh, black] table [z expr=0.0, meta index=2] {data/pb_220/u_approx_1.csv};
    \end{axis}
  \end{tikzpicture} 
   \caption{Solution on $r-$adaptive mesh. \label{fig:test4_sol}}
  \end{subfigure}
  \begin{subfigure}{0.49\textwidth}
	\begin{tikzpicture}  
	\begin{axis}[scale only axis, xlabel = Epoch, ylabel = $\mathcal{L}_{Ritz}(u_{p,r}) - \mathcal{L}_{Ritz}(u_{exact})$, 
	ymode=log,
	ytick pos=left,
	xtick pos=left,
	y label style={at={(-0.17,0.5)}}, 
	height=5.5cm, width=5.5cm, 
	xmin=0, ymin=0.000001, 
	xmax=6000, ymax=0.1,
	legend cell align={left}
	]	
	\addplot[color=black, mark options={solid}] table[x=epoch ,y=lossN]{data/pb_220/loss_clean1.csv};	
	\addplot[color=red, mark options={solid}] table[x=epoch ,y=lossN]{data/pb_220/loss_clean0.csv};	
	\draw[dashed, color=red](1000,0.000001)--(1000,10);
	\node [rotate=90, color = red] at (750,0.0008){\footnotesize \begin{tabular}{c} 1st training step\end{tabular} };	
	\node at (3000,0.02){\footnotesize 2nd training step};
	\end{axis}
	\end{tikzpicture}
   \caption{Loss evolution. \label{fig:test4_loss}}
  \end{subfigure}
\caption{Experiment 5. Approximate solution and loss evolution.} 
\label{fig:Test_4}
\end{figure}
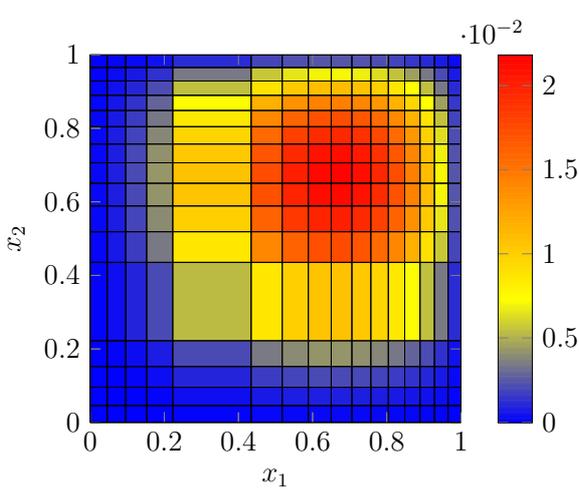
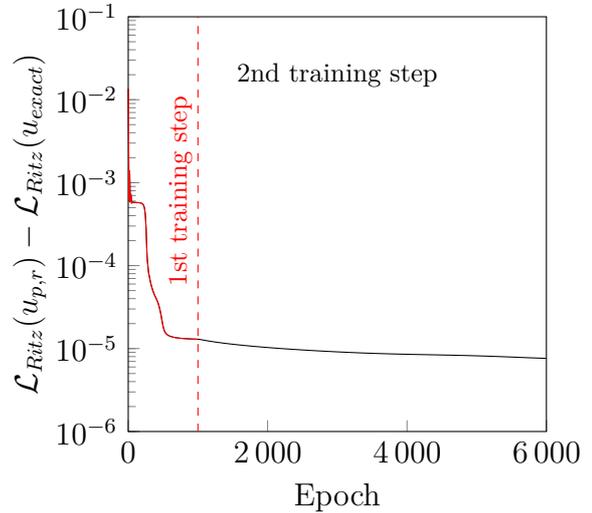

Figure \ref{fig:Test_4} shows the approximated solution on the adapted mesh and the loss evolution for a 2D problem with a smooth solution. 
We compute the results on a 16x16 elements mesh. 
We observe that the $r-$adapted mesh is close to symmetric with respect to the line $y=x$. This, together with the low error in the loss, indicate the high quality of the solution.

\subsection{Experiment 6: L-shape domain}
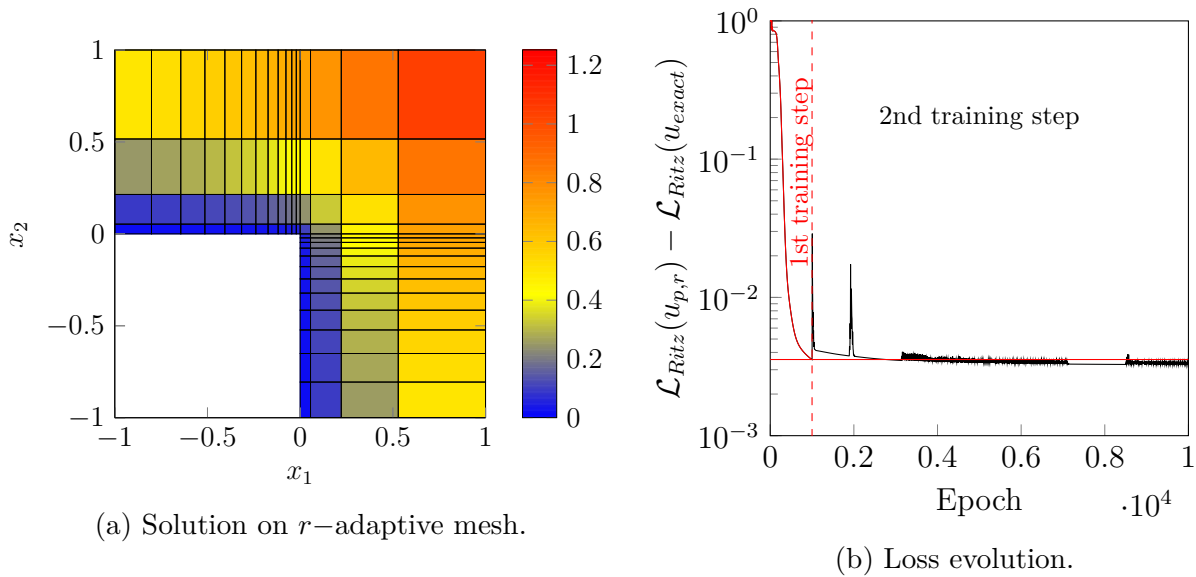
\begin{figure}[!htb]
  \centering
  \begin{subfigure}{0.49\textwidth}
  \begin{tikzpicture}[scale=0.9]
    \begin{axis}[view={0}{90},
    scaled z ticks=false,
	zticklabel=\pgfkeys{/pgf/number format/.cd,fixed,precision=2,zerofill}\pgfmathprintnumber{\tick}, 
    colorbar,
    xlabel = $x_1$,
    ylabel = $x_2$,
    height=7.cm, width=7cm, 
	xmin=-1, ymin=-1, 
	xmax=1, ymax=1,
	zmin = 0, zmax = 0.035,
	ztick = {0,0.01,0.02},
    unbounded coords=jump,
  filter point/.code={%
    \pgfmathparse
      {\pgfkeysvalueof{/data point/x}<0.}%
    \ifpgfmathfloatcomparison
      \pgfmathparse
        {\pgfkeysvalueof{/data point/y}<0.}%
      \ifpgfmathfloatcomparison
        \pgfkeyssetvalue{/data point/x}{nan}%
      \fi
    \fi
  },   
    ]
      \addplot3[surf, point meta=explicit,] table [z expr=0.0, meta index=2] {data/pb_299/u_approx_1.csv};
      \addplot3[surf, point meta=explicit, mesh, black] table [z expr=0.0, meta index=2] {data/pb_299/u_approx_1.csv};
    \end{axis}
  \end{tikzpicture}  
   \caption{Solution on $r-$adaptive mesh. \label{fig:test5_sol}}
  \end{subfigure}
  \begin{subfigure}{0.49\textwidth}
	\begin{tikzpicture}  
	\begin{axis}[scale only axis, xlabel = Epoch, ylabel = $\mathcal{L}_{Ritz}(u_{p,r}) - \mathcal{L}_{Ritz}(u_{exact})$, 
	ymode=log,
	ytick pos=left,
	xtick pos=left,
	y label style={at={(-0.17,0.5)}}, 
	height=5.5cm, width=5.5cm, 
	xmin=0, ymin=0.001, 
	xmax=10000, ymax=1,
	legend cell align={left}
	]	
	\addplot[color=black, mark options={solid}] table[x=epoch ,y=lossN]{data/pb_299/loss_clean1.csv};	
	\addplot[color=red, mark options={solid}] table[x=epoch ,y=lossN]{data/pb_299/loss_clean0.csv};	
	\draw[dashed, color=red](1000,0.000001)--(1000,10);
	\draw[ color=red](0,0.003555)--(10000,0.003555);
	\node [rotate=90, color = red] at (750,0.08){\footnotesize \begin{tabular}{c} 1st training step\end{tabular} };	
	\node at (5000,0.2){\footnotesize 2nd training step};
	\end{axis}
	\end{tikzpicture}
   \caption{Loss evolution. \label{fig:test5_loss}}
  \end{subfigure}
\caption{Experiment 6. Approximate solution and loss evolution.} 
\label{fig:Test_5}
\end{figure}

Figure \ref{fig:Test_5} shows the approximated solution on the adapted mesh and the loss evolution for the so-called L-shape domain problem.

We consider a square reference domain.
To properly define the physical L-shape domain, we fix the points $x_1=0$ and $x_2=0$.
In this way, we ensure the conformity of the elements to the L-shape domain.
The elements outside the domain have null contribution to the loss ---integral equal to zero.
This example shows that we can consider geometries that go beyond simple tensor-product geometries\footnote{Nonetheless, we fall short to consider arbitrary geometries}.
Figure \ref{fig:test5_sol} shows that the final mesh computed over the physical domain is almost symmetric with respect to the line $x=y$ and has fewer elements than the uniform initial one. 
The reference square domain employs two 16-dimensional vectors of 1D coordinates (producing 192 elements). In contrast, the final adapted mesh only employs 112 elements in the L-shape area of interest (the rest are zero). 
However, the solution obtained in the adapted mesh constrains lower error than the initial one, as observed by the reduction in the loss value. In the $r-$adapted mesh, the element sizes grow exponentially on both sides from the singularity, increasing the density of elements in the null area of the domain (outside the L-shape).

\section{Conclusions and Future Work}
We propose a DL $r-$adaptive method to solve PDEs.
The method simultaneously optimizes the mesh and approximated solution.
We have implemented the method for one and two spatial dimensions using tensor product meshes.
Numerical experiments show promising results for solutions that are smooth, singular, or exhibit strong gradients. 
The method outperforms the use of uniform meshes.

In this work, we have considered piecewise-linear functions to approximate the solution. The extension to higher-degree piecewise-polynomial functions is straightforward.
Moreover, we have restricted meshes with tensor-product topology.
The extension to triangular/tetrahedral meshes will be pursued as part of our future effort in the area.

\section*{Acknowledgments}
David Pardo has received funding from: 
the European Union's Horizon 2020 research and innovation program under the Marie Sklodowska-Curie grant agreement No 777778 (MATHROCKS); 
the Spanish Ministry of Science and Innovation projects with references TED2021-132783B-I00, PID2019-108111RB-I00 (FEDER/AEI) and PDC2021-121093-I00 (AEI/Next Generation EU), 
the ``BCAM Severo Ochoa'' accreditation of excellence (SEV-2017-0718); 
and the Basque Government through the BERC 2022-2025 program, 
the three Elkartek projects 3KIA (KK-2020/00049), EXPERTIA (KK-2021/00048), and SIGZE (KK-2021/00095), 
and the Consolidated Research Group MATHMODE (IT1456-22) given by the Department of Education.

\bibliographystyle{elsarticle-num}
\bibliography{biblio}

\begin{thebibliography}{10}
\expandafter\ifx\csname url\endcsname\relax
  \def\url#1{\texttt{#1}}\fi
\expandafter\ifx\csname urlprefix\endcsname\relax\def\urlprefix{URL }\fi
\expandafter\ifx\csname href\endcsname\relax
  \def\href#1#2{#2} \def\path#1{#1}\fi

\bibitem{bookDL}
I.~Goodfellow, Y.~Bengio, A.~Courville, Deep Learning, MIT Press, 2016,
  \url{http://www.deeplearningbook.org}.

\bibitem{Koll2021}
S.~Kollmannsberger, D.~D'Angella, M.~Jokeit, L.~Herrmann,
  \href{https://doi.org/10.1007/978-3-030-76587-3}{Deep Learning in
  Computational Mechanics: An itroductory Course}, Springer Cham, 2021.
\newblock \href {https://doi.org/10.1007/978-3-030-76587-3}
  {\path{doi:10.1007/978-3-030-76587-3}}.
\newline\urlprefix\url{https://doi.org/10.1007/978-3-030-76587-3}

\bibitem{DLapplications}
L.~Deng, D.~Yu, \href{https://doi.org/10.1561/2000000039}{Deep {L}earning:
  {M}ethods and {A}pplications}, Found. Trends Signal Process. 7~(3–4) (2014)
  197–387.
\newblock \href {https://doi.org/10.1561/2000000039}
  {\path{doi:10.1561/2000000039}}.
\newline\urlprefix\url{https://doi.org/10.1561/2000000039}

\bibitem{biomedicine}
M.~Wainberg, D.~Merico, A.~Delong, B.~J. Frey, Deep learning in biomedicine,
  Nature biotechnology 36~(9) (2018) 829--838.
\newblock \href {https://doi.org/https://doi.org/10.1038/nbt.4233}
  {\path{doi:https://doi.org/10.1038/nbt.4233}}.

\bibitem{puente_ana}
A.~Fernandez-Navamuel, F.~Magalhães, D.~Zamora-Sánchez, Ángel J~Omella,
  D.~Garcia-Sanchez, D.~Pardo,
  \href{https://doi.org/10.1177/14759217211041684}{Deep learning enhanced
  principal component analysis for structural health monitoring}, Structural
  Health Monitoring 0~(0) (2022) 14759217211041684.
\newblock \href
  {http://arxiv.org/abs/https://doi.org/10.1177/14759217211041684}
  {\path{arXiv:https://doi.org/10.1177/14759217211041684}}, \href
  {https://doi.org/10.1177/14759217211041684}
  {\path{doi:10.1177/14759217211041684}}.
\newline\urlprefix\url{https://doi.org/10.1177/14759217211041684}

\bibitem{Ana_puentes_FEM}
A.~Fernandez-Navamuel, D.~Zamora-Sánchez, Ángel J.~Omella, D.~Pardo,
  D.~Garcia-Sanchez, F.~Magalhães,
  \href{https://www.sciencedirect.com/science/article/pii/S0141029622001638}{Supervised
  deep learning with finite element simulations for damage identification in
  bridges}, Engineering Structures 257 (2022) 114016.
\newblock \href
  {https://doi.org/https://doi.org/10.1016/j.engstruct.2022.114016}
  {\path{doi:https://doi.org/10.1016/j.engstruct.2022.114016}}.
\newline\urlprefix\url{https://www.sciencedirect.com/science/article/pii/S0141029622001638}

\bibitem{extra_deep}
S.~Alyaev, M.~Shahriari, D.~Pardo, Ángel Javier~Omella, D.~S. Larsen,
  N.~Jahani, E.~Suter, \href{https://doi.org/10.1190/geo2020-0389.1}{Modeling
  extra-deep electromagnetic logs using a deep neural network}, GEOPHYSICS
  86~(3) (2021) E269--E281.
\newblock \href {http://arxiv.org/abs/https://doi.org/10.1190/geo2020-0389.1}
  {\path{arXiv:https://doi.org/10.1190/geo2020-0389.1}}, \href
  {https://doi.org/10.1190/geo2020-0389.1} {\path{doi:10.1190/geo2020-0389.1}}.
\newline\urlprefix\url{https://doi.org/10.1190/geo2020-0389.1}

\bibitem{BREVIS}
I.~Brevis, I.~Muga, K.~G. {van der Zee},
  \href{https://www.sciencedirect.com/science/article/pii/S0898122120303199}{A
  machine-learning minimal-residual ({ML}-{MR}es) framework for goal-oriented
  finite element discretizations}, Computers \& Mathematics with Applications
  95 (2021) 186--199, recent Advances in Least-Squares and Discontinuous
  Petrov–Galerkin Finite Element Methods.
\newblock \href {https://doi.org/https://doi.org/10.1016/j.camwa.2020.08.012}
  {\path{doi:https://doi.org/10.1016/j.camwa.2020.08.012}}.
\newline\urlprefix\url{https://www.sciencedirect.com/science/article/pii/S0898122120303199}

\bibitem{Brevis2}
I.~Brevis, I.~Muga, K.~G. van~der Zee,
  \href{https://arxiv.org/abs/2206.07475}{Neural control of discrete weak
  formulations: Galerkin, least-squares and minimal-residual methods with
  quasi-optimal weights} (2022).
\newblock \href {https://doi.org/10.48550/ARXIV.2206.07475}
  {\path{doi:10.48550/ARXIV.2206.07475}}.
\newline\urlprefix\url{https://arxiv.org/abs/2206.07475}

\bibitem{Maciej1}
M.~Paszynski, R.~Grzeszczuk, D.~Pardo, L.~Demkowicz, Deep learning driven
  self-adaptive hp finite element method, in: M.~Paszynski,
  D.~Kranzlm{\"u}ller, V.~V. Krzhizhanovskaya, J.~J. Dongarra, P.~M.~A. Sloot
  (Eds.), Computational Science -- ICCS 2021, Springer International
  Publishing, Cham, 2021, pp. 114--121.

\bibitem{Maciej2}
T.~Sluzalec, R.~Grzeszczuk, S.~Rojas, W.~Dzwinel, M.~Paszynski,
  \href{https://arxiv.org/abs/2209.05844}{Quasi-optimal $hp$-finite element
  refinements towards singularities via deep neural network prediction} (2022).
\newblock \href {https://doi.org/10.48550/ARXIV.2209.05844}
  {\path{doi:10.48550/ARXIV.2209.05844}}.
\newline\urlprefix\url{https://arxiv.org/abs/2209.05844}

\bibitem{carlos}
C.~Uriarte, D.~Pardo, Ángel Javier~Omella,
  \href{https://www.sciencedirect.com/science/article/pii/S0045782521007374}{A
  {F}inite {E}lement based {D}eep {L}earning solver for parametric {PDEs}},
  Computer Methods in Applied Mechanics and Engineering 391 (2022) 114562.
\newblock \href {https://doi.org/https://doi.org/10.1016/j.cma.2021.114562}
  {\path{doi:https://doi.org/10.1016/j.cma.2021.114562}}.
\newline\urlprefix\url{https://www.sciencedirect.com/science/article/pii/S0045782521007374}

\bibitem{overview}
C.~Beck, M.~Hutzenthaler, A.~Jentzen, B.~Kuckuck, An overview on deep
  learning-based approximation methods for partial differential equations
  (2021).
\newblock \href {http://arxiv.org/abs/2012.12348} {\path{arXiv:2012.12348}}.

\bibitem{PIML}
G.~Karniadakis, I.~Kevrekidis, L.~Lu, P.~Perdikaris, S.~Wang, L.Yang,
  Physics-informed machine learning, Nature Reviews Physics 3~(6) (2021)
  422--440.
\newblock \href {https://doi.org/https://doi.org/10.1038/s42254-021-00314-5}
  {\path{doi:https://doi.org/10.1038/s42254-021-00314-5}}.

\bibitem{PINN}
M.~Raissi, P.~Perdikaris, G.~Karniadakis,
  \href{https://www.sciencedirect.com/science/article/pii/S0021999118307125}{Physics-informed
  neural networks: A deep learning framework for solving forward and inverse
  problems involving nonlinear partial differential equations}, Journal of
  Computational Physics 378 (2019) 686--707.
\newblock \href {https://doi.org/https://doi.org/10.1016/j.jcp.2018.10.045}
  {\path{doi:https://doi.org/10.1016/j.jcp.2018.10.045}}.
\newline\urlprefix\url{https://www.sciencedirect.com/science/article/pii/S0021999118307125}

\bibitem{DeepRitz}
E.~Weinan, Y.~Bing, \href{https://doi.org/10.1007/s40304-018-0127-z}{The deep
  {R}itz {M}ethod: {A} {D}eep {L}earning-{B}ased {N}umerical {A}lgorithm for
  {S}olving {V}ariational {P}roblems}, Communications in Mathematics and
  Statistics 6 (2018) 1--12.
\newblock \href {https://doi.org/10.1007/s40304-018-0127-z}
  {\path{doi:10.1007/s40304-018-0127-z}}.
\newline\urlprefix\url{https://doi.org/10.1007/s40304-018-0127-z}

\bibitem{DeepLS}
Z.~Cai, J.~Chen, M.~Liu, X.~Liu,
  \href{https://www.sciencedirect.com/science/article/pii/S0021999120304812}{Deep
  least-squares methods: An unsupervised learning-based numerical method for
  solving elliptic {PDEs}}, Journal of Computational Physics 420 (2020) 109707.
\newblock \href {https://doi.org/https://doi.org/10.1016/j.jcp.2020.109707}
  {\path{doi:https://doi.org/10.1016/j.jcp.2020.109707}}.
\newline\urlprefix\url{https://www.sciencedirect.com/science/article/pii/S0021999120304812}

\bibitem{DGM}
J.~Sirignano, K.~Spiliopoulos,
  \href{https://www.sciencedirect.com/science/article/pii/S0021999118305527}{{DGM}:
  A deep learning algorithm for solving partial differential equations},
  Journal of Computational Physics 375 (2018) 1339--1364.
\newblock \href {https://doi.org/https://doi.org/10.1016/j.jcp.2018.08.029}
  {\path{doi:https://doi.org/10.1016/j.jcp.2018.08.029}}.
\newline\urlprefix\url{https://www.sciencedirect.com/science/article/pii/S0021999118305527}

\bibitem{hpVPINNs}
E.~Kharazmi, Z.~Zhang, G.~E. Karniadakis,
  \href{https://www.sciencedirect.com/science/article/pii/S0045782520307325}{hp-{VPINN}s:
  Variational physics-informed neural networks with domain decomposition},
  Computer Methods in Applied Mechanics and Engineering 374 (2021) 113547.
\newblock \href {https://doi.org/https://doi.org/10.1016/j.cma.2020.113547}
  {\path{doi:https://doi.org/10.1016/j.cma.2020.113547}}.
\newline\urlprefix\url{https://www.sciencedirect.com/science/article/pii/S0045782520307325}

\bibitem{high_dimensions2}
T.~Poggio, H.~Mhaskar, L.~Rosasco, B.~Miranda, Q.~Liao,
  \href{https://doi.org/10.1007/s11633-017-1054-2}{Why and when can deep-but
  not shallow-networks avoid the curse of dimensionality: A review},
  International Journal of Automation and Computing 14~(5) (2017) 503--519.
\newblock \href {https://doi.org/10.1007/s11633-017-1054-2}
  {\path{doi:10.1007/s11633-017-1054-2}}.
\newline\urlprefix\url{https://doi.org/10.1007/s11633-017-1054-2}

\bibitem{high_dimensions}
J.~Han, A.~Jentzen, W.~E,
  \href{https://www.pnas.org/content/115/34/8505}{Solving high-dimensional
  partial differential equations using deep learning}, Proceedings of the
  National Academy of Sciences 115~(34) (2018) 8505--8510.
\newblock \href
  {http://arxiv.org/abs/https://www.pnas.org/content/115/34/8505.full.pdf}
  {\path{arXiv:https://www.pnas.org/content/115/34/8505.full.pdf}}, \href
  {https://doi.org/10.1073/pnas.1718942115}
  {\path{doi:10.1073/pnas.1718942115}}.
\newline\urlprefix\url{https://www.pnas.org/content/115/34/8505}

\bibitem{Montecarlo_no_lineal}
W.~E, J.~Han, A.~Jentzen,
  \href{https://doi.org/10.1088/1361-6544/ac337f}{Algorithms for solving high
  dimensional {PDEs}: from nonlinear {M}onte {C}arlo to machine learning},
  Nonlinearity 35~(1) (2021) 278--310.
\newblock \href {https://doi.org/10.1088/1361-6544/ac337f}
  {\path{doi:10.1088/1361-6544/ac337f}}.
\newline\urlprefix\url{https://doi.org/10.1088/1361-6544/ac337f}

\bibitem{moin_2010}
P.~Moin, Fundamentals of Engineering Numerical Analysis, 2nd Edition, Cambridge
  University Press, 2010.
\newblock \href {https://doi.org/10.1017/CBO9780511781438}
  {\path{doi:10.1017/CBO9780511781438}}.

\bibitem{JANDER_integration}
J.~A. Rivera, J.~M. Taylor, Ángel J.~Omella, D.~Pardo,
  \href{https://www.sciencedirect.com/science/article/pii/S0045782522000810}{On
  quadrature rules for solving partial differential equations using neural
  networks}, Computer Methods in Applied Mechanics and Engineering 393 (2022)
  114710.
\newblock \href {https://doi.org/https://doi.org/10.1016/j.cma.2022.114710}
  {\path{doi:https://doi.org/10.1016/j.cma.2022.114710}}.
\newline\urlprefix\url{https://www.sciencedirect.com/science/article/pii/S0045782522000810}

\bibitem{adaptatividad_origienes}
W.~Carroll, R.~Barker,
  \href{https://www.sciencedirect.com/science/article/pii/0020768373900115}{A
  theorem for optimum finite-element idealizations}, International Journal of
  Solids and Structures 9~(7) (1973) 883--895.
\newblock \href {https://doi.org/https://doi.org/10.1016/0020-7683(73)90011-5}
  {\path{doi:https://doi.org/10.1016/0020-7683(73)90011-5}}.
\newline\urlprefix\url{https://www.sciencedirect.com/science/article/pii/0020768373900115}

\bibitem{adaptatividad_origienes1}
H.~G.~W. Burchard, D.~F. Hale, Piecewise polynomial approximation on optimal
  meshes, Journal of Approximation Theory 14 (1975) 128--147.

\bibitem{adaptatividad_origienes2}
I.~Babu{\v{s}}ka, W.~C. Rheinboldt,
  \href{http://inis.iaea.org/search/search.aspx?orig_q=RN:10492324}{Analysis of
  optimal finite-element meshes in {R$^1$}}, Mathematics of Computation
  33~(146) (1979) 435--463.
\newline\urlprefix\url{http://inis.iaea.org/search/search.aspx?orig_q=RN:10492324}

\bibitem{Alfonzetti}
S.~Alfonzetti, A finite element mesh generator based on an adaptive neural
  network, IEEE Transactions on Magnetics 34~(5) (1998) 3363--3366.
\newblock \href {https://doi.org/10.1109/20.717791}
  {\path{doi:10.1109/20.717791}}.

\bibitem{MANEVITZ2005447}
L.~Manevitz, A.~Bitar, D.~Givoli,
  \href{https://www.sciencedirect.com/science/article/pii/S0925231204003078}{Neural
  network time series forecasting of finite-element mesh adaptation},
  Neurocomputing 63 (2005) 447--463, new Aspects in Neurocomputing: 11th
  European Symposium on Artificial Neural Networks.
\newblock \href {https://doi.org/https://doi.org/10.1016/j.neucom.2004.06.009}
  {\path{doi:https://doi.org/10.1016/j.neucom.2004.06.009}}.
\newline\urlprefix\url{https://www.sciencedirect.com/science/article/pii/S0925231204003078}

\bibitem{BOHN202161}
J.~Bohn, M.~Feischl,
  \href{https://www.sciencedirect.com/science/article/pii/S089812212100198X}{Recurrent
  neural networks as optimal mesh refinement strategies}, Computers \&
  Mathematics with Applications 97 (2021) 61--76.
\newblock \href {https://doi.org/https://doi.org/10.1016/j.camwa.2021.05.018}
  {\path{doi:https://doi.org/10.1016/j.camwa.2021.05.018}}.
\newline\urlprefix\url{https://www.sciencedirect.com/science/article/pii/S089812212100198X}

\bibitem{budd_huang_russell_2009}
C.~J. Budd, W.~Huang, R.~D. Russell, Adaptivity with moving grids, Acta
  Numerica 18 (2009) 111–241.
\newblock \href {https://doi.org/10.1017/S0962492906400015}
  {\path{doi:10.1017/S0962492906400015}}.

\bibitem{DORFI1987175}
E.~Dorfi, L.~Drury,
  \href{https://www.sciencedirect.com/science/article/pii/0021999187901616}{Simple
  adaptive grids for 1 - d initial value problems}, Journal of Computational
  Physics 69~(1) (1987) 175--195.
\newblock \href {https://doi.org/https://doi.org/10.1016/0021-9991(87)90161-6}
  {\path{doi:https://doi.org/10.1016/0021-9991(87)90161-6}}.
\newline\urlprefix\url{https://www.sciencedirect.com/science/article/pii/0021999187901616}

\bibitem{huang2011adaptive}
W.~Huang, R.~Russell, Adaptive {M}oving {M}esh {M}ethods, Applied Mathematical
  Sciences, Springer, New York, NY, 2011.
\newblock \href {https://doi.org/https://doi.org/10.1007/978-1-4419-7916-2}
  {\path{doi:https://doi.org/10.1007/978-1-4419-7916-2}}.

\bibitem{MMPDE}
W.~Huang, Y.~Ren, R.~D. Russell, \href{https://doi.org/10.1137/0731038}{Moving
  {M}esh {P}artial {D}ifferential {E}quations ({MMPDES}) {B}ased on the
  {E}quidistribution {P}rinciple}, SIAM Journal on Numerical Analysis 31~(3)
  (1994) 709--730.
\newblock \href {http://arxiv.org/abs/https://doi.org/10.1137/0731038}
  {\path{arXiv:https://doi.org/10.1137/0731038}}, \href
  {https://doi.org/10.1137/0731038} {\path{doi:10.1137/0731038}}.
\newline\urlprefix\url{https://doi.org/10.1137/0731038}

\bibitem{Budd2}
C.~J. Budd, A.~T. McRae, C.~J. Cotter,
  \href{https://www.sciencedirect.com/science/article/pii/S0021999118305515}{The
  scaling and skewness of optimally transported meshes on the sphere}, Journal
  of Computational Physics 375 (2018) 540--564.
\newblock \href {https://doi.org/https://doi.org/10.1016/j.jcp.2018.08.028}
  {\path{doi:https://doi.org/10.1016/j.jcp.2018.08.028}}.
\newline\urlprefix\url{https://www.sciencedirect.com/science/article/pii/S0021999118305515}

\bibitem{Autodiff}
B.~v. Merri\"{e}nboer, O.~Breuleux, A.~Bergeron, P.~Lamblin, Automatic
  {D}ifferentiation in {ML}: {W}here {W}e {A}re and {W}here {W}e {S}hould {B}e
  {G}oing, in: Proceedings of the 32nd International Conference on Neural
  Information Processing Systems, NIPS'18, Curran Associates Inc., Red Hook,
  NY, USA, 2018, p. 8771–8781.

\bibitem{Autodiff_1}
A.~G. Baydin, B.~A. Pearlmutter, A.~A. Radul, J.~M. Siskind, Automatic
  {D}ifferentiation in {M}achine {L}earning: {A} {S}urvey, J. Mach. Learn. Res.
  18~(1) (2017) 5595–5637.

\bibitem{Autodiff_2}
A.~Griewank, A.~Walther, \href{http://bookstore.siam.org/ot105/}{Evaluating
  Derivatives: {P}rinciples and Techniques of Algorithmic Differentiation}, 2nd
  Edition, no. 105 in Other Titles in Applied Mathematics, SIAM, Philadelphia,
  PA, 2008.
\newline\urlprefix\url{http://bookstore.siam.org/ot105/}

\bibitem{backprop}
D.~E. Rumelhart, G.~E. Hinton, R.~J. Williams, Learning representations by
  back-propagating errors, Nature 323 (1986) 533--536.
\newblock \href {https://doi.org/https://doi.org/10.1038/323533a0}
  {\path{doi:https://doi.org/10.1038/323533a0}}.

\bibitem{mesh}
P.~J. Frey, P.~L. George,
  \href{https://onlinelibrary.wiley.com/doi/abs/10.1002/9780470611166}{Mesh
  Generation}, John Wiley \& Sons, Ltd, 2008.
\newblock \href {https://doi.org/https://doi.org/10.1002/9780470611166}
  {\path{doi:https://doi.org/10.1002/9780470611166}}.
\newline\urlprefix\url{https://onlinelibrary.wiley.com/doi/abs/10.1002/9780470611166}

\bibitem{cell_method}
J.~Parvizian, A.~D{\"u}ster, E.~Rank,
  \href{https://doi.org/10.1007/s00466-007-0173-y}{Finite cell method},
  Computational Mechanics 41~(1) (2007) 121--133.
\newblock \href {https://doi.org/10.1007/s00466-007-0173-y}
  {\path{doi:10.1007/s00466-007-0173-y}}.
\newline\urlprefix\url{https://doi.org/10.1007/s00466-007-0173-y}

\bibitem{cell_method_2}
N.~Zander, S.~Kollmannsberger, M.~Ruess, Z.~Yosibash, E.~Rank,
  \href{https://www.sciencedirect.com/science/article/pii/S0898122112005688}{The
  finite cell method for linear thermoelasticity}, Computers \& Mathematics
  with Applications 64~(11) (2012) 3527--3541.
\newblock \href {https://doi.org/https://doi.org/10.1016/j.camwa.2012.09.002}
  {\path{doi:https://doi.org/10.1016/j.camwa.2012.09.002}}.
\newline\urlprefix\url{https://www.sciencedirect.com/science/article/pii/S0898122112005688}

\bibitem{PaulChew1989}
L.~Paul~Chew, \href{https://doi.org/10.1007/BF01553881}{Constrained delaunay
  triangulations}, Algorithmica 4~(1) (1989) 97--108.
\newblock \href {https://doi.org/10.1007/BF01553881}
  {\path{doi:10.1007/BF01553881}}.
\newline\urlprefix\url{https://doi.org/10.1007/BF01553881}

\bibitem{tensorflow2015-whitepaper}
M.~Abadi, A.~Agarwal, P.~Barham, E.~Brevdo, Z.~Chen, C.~Citro, G.~S. Corrado,
  A.~Davis, J.~Dean, M.~Devin, S.~Ghemawat, I.~Goodfellow, A.~Harp, G.~Irving,
  M.~Isard, Y.~Jia, R.~Jozefowicz, L.~Kaiser, M.~Kudlur, J.~Levenberg,
  D.~Man\'{e}, R.~Monga, S.~Moore, D.~Murray, C.~Olah, M.~Schuster, J.~Shlens,
  B.~Steiner, I.~Sutskever, K.~Talwar, P.~Tucker, V.~Vanhoucke, V.~Vasudevan,
  F.~Vi\'{e}gas, O.~Vinyals, P.~Warden, M.~Wattenberg, M.~Wicke, Y.~Yu,
  X.~Zheng, \href{https://www.tensorflow.org/}{{TensorFlow}: Large-scale
  machine learning on heterogeneous systems}, software available from
  tensorflow.org (2015).
\newline\urlprefix\url{https://www.tensorflow.org/}

\bibitem{maclaurin2015autograd}
D.~Maclaurin, D.~Duvenaud, R.~P. Adams, Autograd: Effortless gradients in
  numpy, in: ICML 2015 AutoML Workshop, Vol. 238, 2015, p.~5.

\bibitem{jax_github}
J.~Bradbury, R.~Frostig, P.~Hawkins, M.~J. Johnson, C.~Leary, D.~Maclaurin,
  G.~Necula, A.~Paszke, J.~Vander{P}las, S.~Wanderman-{M}ilne, Q.~Zhang,
  \href{http://github.com/google/jax}{{JAX}: composable transformations of
  {P}ython+{N}um{P}y programs} (2022).
\newline\urlprefix\url{http://github.com/google/jax}

\bibitem{Ritz}
W.~Ritz, \href{http://eudml.org/doc/149295}{\"{U}ber eine neue methode zur
  l\"osung gewisser variationsprobleme der mathematischen physik.}, Journal
  f\"ur die reine und angewandte Mathematik 135 (1909) 1--61.
\newline\urlprefix\url{http://eudml.org/doc/149295}

\bibitem{Claes}
C.~Johnson, Numerical {S}olution of {P}artial {D}ifferential {E}quations by the
  {F}inite {E}lement {M}ethod, Cambridge U. Press, 1987.

\bibitem{adam}
D.~P. Kingma, J.~Ba, \href{https://arxiv.org/abs/1412.6980}{Adam: A method for
  stochastic optimization} (2014).
\newblock \href {https://doi.org/10.48550/ARXIV.1412.6980}
  {\path{doi:10.48550/ARXIV.1412.6980}}.
\newline\urlprefix\url{https://arxiv.org/abs/1412.6980}

\bibitem{Strouboulis}
I.~Babu{\v{s}}ka, S.~Theofanis, {T}he {F}inite {E}lement {M}ethod and {I}ts
  {R}eliability, Numerical mathematics and scientific computation, New York:
  Oxford University Press, 2001.

\bibitem{Gui1986b}
W.~Gui, I.~Babu{\v{s}}ka, The $h$, $p$ and $h-p$ versions of the finite element
  method in 1 dimension: Part {II}. the error analysis of the $h-$ and $h-p$
  versions, Numerische Mathematik 49~(6) (1986) 613--657.

\end{thebibliography}
\end{document}